\documentclass[twoside,11pt,reqno]{amsart}
\usepackage{amsmath,amssymb,amscd,mathrsfs}
\usepackage{graphics,verbatim}

\newcommand{\p}[1]{\ensuremath{\overline{#1}}}
\newcommand{\losemi}{{\otimes \kern -.78em \ltimes}}
\newcommand{\rosemi}{{\otimes \kern -.78em \rtimes}}

\newcommand{\Hom}{\ensuremath{\operatorname{Hom}}}

\newcommand{\Ind}{\ensuremath{\operatorname{ind}}}

\newcommand{\Ker}{\ensuremath{\operatorname{Ker} }}

\newcommand{\sgn}{\operatorname{sgn}}
\newcommand{\Ext}{\operatorname{Ext}}

\newcommand{\0}{\bar 0}
\newcommand{\1}{\bar 1}
\newcommand{\Z}{\mathbb{Z}}
\newcommand{\C}{\mathbb{C}}

\newcommand{\setof}[2]{\ensuremath{\left\{ #1 \:|\: #2 \right\}}}
\newcommand{\gl}{\ensuremath{\mathfrak{gl}(m|n)}}
\newcommand{\g}{\ensuremath{\mathfrak{g}}}
\newcommand{\e}{\ensuremath{\mathfrak{f}}}
\renewcommand{\c}{\ensuremath{\mathfrak{e}}}
\newcommand{\x}{\ensuremath{\langle x \rangle}}

\newcommand{\res}{\ensuremath{\operatorname{res}}}

\newcommand{\fg}{\ensuremath{\mathfrak{g}}}
\newcommand{\fh}{\ensuremath{\mathfrak{h}}}
\newcommand{\fc}{\ensuremath{\c}}
\newcommand{\fe}{\ensuremath{\e}}
\newcommand{\fq}{\ensuremath{\mathfrak{q}}}
\renewcommand{\a}{\alpha}
\renewcommand{\b}{\beta}
\renewcommand{\d}{\delta}
\newcommand{\ep}{\epsilon}
\newcommand{\ve}{\varepsilon}
\newcommand{\ga}{\gamma}
\newcommand{\s}{\sigma}
\newcommand{\slmn}{\ensuremath{\mathfrak{sl}(m|n)}}
\newcommand{\slnn}{\ensuremath{\mathfrak{psl}(n|n)}}
\newcommand{\osp}{\ensuremath{\mathfrak{osp}}}
\newcommand{\str}{\ensuremath{\operatorname{str} }}
\newcommand{\tr}{\ensuremath{\operatorname{tr} }}

\newcommand{\iz}{\ensuremath{i_{0}}}

\def\Label{\label}

\newtheorem{Df}{Definition}[subsection]
\newtheorem{theorem}[Df]{Theorem}

\newtheorem{lemma}[Df]{Lemma}

\newtheorem{corollary}[Df]{Corollary}
\newtheorem{prop}[Df]{Proposition}
\newtheorem{conj}[Df]{Conjecture} 

\numberwithin{equation}{subsection}

\def\lierelcohom{\ref{SS:lierelcohom}}
\def\Lrelatescohom2{\ref{L:relatescohom2}}

\def\CORextequivalence{\ref{C:extequivalence}}

\def\THcohomring{\ref{T:cohomring}}
\def\THfinitegen2{\ref{T:finitegen2}} 
 
\def\SSLR{\ref{SS:stableactions}}
\def\THlunarichardson{\ref{T:lunarichardson}}

\def\THdadokkac{\ref{T:dadokkac}}
 
\def\THKosRal{\ref{T:KR}}
\def\THKRforg{\ref{T:KRforg}}
\def\THKRfore{\ref{T:KRfore}}
\def\THKRforc{\ref{T:KRforc}}

\def\THinvertingD{\ref{T:invertingD}}
\def\Smorecohom{\ref{S:morecohom}}

\def\THinjectivityoncohains{\ref{T:injectivityoncohains}}

\def\THkrulldimoneiso{\ref{T:krulldimoneiso}}
\def\THlocalizingcohom{\ref{T:localizingcohom}}

\def\LirrepsofQ(1){\ref{L:irrepsofQ(1)}}

\def\Crankcor1{\ref{C:rankcor1} }

\begin{document}
\title{Cohomology and Support Varieties for Lie Superalgebras}

\author{Brian D. Boe }
\address{Department of Mathematics \\
            University of Georgia \\
            Athens, GA 30602}
\email{brian@math.uga.edu, kujawa@math.uga.edu, nakano@math.uga.edu}
\thanks{Research of the first author was partially supported by NSA 
grant H98230-04-1-0103}
\author{Jonathan R. Kujawa}
\thanks{Research of the second author was partially supported by NSF grant
DMS-0402916}\
\author{Daniel K. Nakano}
\thanks{Research of the third author was partially supported by NSF
grant  DMS-0400548}
\date{\today}
%\subjclass[2000]{Primary 17B56, 14E20; Secondary 46E25, 20C20}

\maketitle

\section{Introduction} \Label{S:intro}

\subsection{}\setcounter{Df}{0} The blocks of the Category ${\mathcal O}$ (or relative 
Category ${\mathcal O}_{S}$) for complex semisimple Lie algebras are 
well known examples of highest weight categories, as defined in \cite{CPS}, 
with finitely many simple modules. These facts imply that the projective 
resolutions for modules in these categories have finite length, so the cohomology (or 
extensions) can be non-zero in only finitely many degrees. 

On the other hand, one can consider the category $\mathcal{F}$ of finite-dimensional supermodules for 
a classical Lie superalgebra $\g$ (e.g.\ ${\mathfrak g}=\mathfrak{gl}(m|n)$) which are completely 
reducible over ${\mathfrak g}_{\0}$. For $\mathfrak{gl}(m|n)$ this was shown to be a highest weight category 
by Brundan \cite{brundan} and recently also for $\mathfrak{osp}(2|2n)$ by Cheng, Wang, and Zhang \cite{CWZ}. 
However, in these categories there are infinitely many simple modules and   
projective resolutions often have infinite length.  Cohomology can be non-zero in 
infinitely many degrees and in fact can have polynomial rate of growth. 

A similar scenario 
has been successfully handled for modular representations of finite groups and restricted Lie algebras by 
the use of support varieties.   Remarkably, an analogous theory holds for $\g$-supermodules in $\mathcal{F}$.  In this paper we develop this algebro-geometric approach and explore the connections it provides between the representation theory, cohomology, and combinatorics of classical Lie superalgebras and the category $\mathcal{F}.$  

\subsection{}\setcounter{Df}{0}
Our approach incorporates many previously 
known computations related to relative cohomology for Lie superalgebras by Brundan \cite{brundan}, 
Fuks and Leites \cite{leitesfuks}, and Gruson \cite{Gru:97,Gru:00}.  It should be noted 
that it is important to use relative cohomology for our purposes rather than
ordinary Lie superalgebra cohomology.  The latter is usually non-zero in only   
finitely many degrees (cf.\ \cite{fuks,Gru:97}) and thus will not capture much information about 
the representation theory of ${\mathfrak g}$.  

Let us set some notation.  Let ${\mathfrak g}={\mathfrak g}_{\0}\oplus {\mathfrak g}_{\1}$ be a 
simple classical Lie superalgebra over the complex numbers as classified by Kac \cite{K}.  A description of the simple classical Lie superalgebras can be found in Sections~\ref{SS:typeA}-\ref{SS:typeF4}.   Denote by $G_{\0}$ the connected 
reductive algebraic group such that $\operatorname{Lie}(G_{\0})={\mathfrak g}_{\0}$.  

Perhaps the most striking outcome of our investigation is that, under certain mild assumptions, $\g$ contains Lie subsuperalgebras which ``detect'' the cohomology of $\g.$   In Section~\ref{S:morecohom} we construct these Lie subsuperalgebras 
of ${\mathfrak g}$ and show that they arise naturally from results in the 
invariant theory of reductive groups due to Luna and Richardson \cite{luna}, and 
Dadok and Kac \cite{dadokkac}. In particular, if $R=\operatorname{H}^{\bullet}
({\mathfrak g},{\mathfrak g}_{\0};{\mathbb  C})$ is 
the cohomology for ${\mathfrak g}$ relative 
to ${\mathfrak g}_{\0},$ then there exists a Lie subsuperalgebra 
${\mathfrak e}={\mathfrak e}_{\0}\oplus {\mathfrak e}_{\1}$ such that 
\begin{equation}\label{E:introiso}
R\cong S^{\bullet}({\mathfrak g}_{\1}^{*})^{G_{\0}}\cong S^{\bullet}({\mathfrak e}_{\1}^{*})^{W} 
\cong \text{H}^{\bullet}({\mathfrak e},{\mathfrak e}_{\0};{\mathbb C})^{W}
\end{equation} 
where $W$ is a finite pseudoreflection group. From this isomorphism one sees that 
$R$ is a finitely generated algebra and in fact a polynomial algebra.  To see an example of the subsuperalgebra $\c$, we refer the reader to Sections~\ref{SS:DK}-\ref{SS:LR} where the type A case is considered in some detail.

The subsuperalgebra $\c$ can be viewed 
as an analogue of a Sylow subgroup with the above theorem looking very much like a 
theorem involving transfer for finite groups.  The 
reader may also notice that the aforementioned isomorphism is reminiscent 
of Borel's calculation of the rational cohomology for the classifying space 
of a compact Lie group.

In recent years there has been active and ongoing interest in 
the development and use of support variety theories in various contexts  
(cf.\ Balmer \cite{balmer}, Erdmann-Holloway \cite{erdmannholloway}, 
Friedlander-Pevtsova \cite{friedlanderpevtsova}, Snashall-Solberg \cite{snashall}). By using the finite generation of $R$ we develop such a theory for ${\mathfrak g}$-supermodules in the category $\mathcal{F}$. 
Given such a ${\mathfrak g}$-supermodule $M$ one can consider the support 
varieties ${\mathcal V}_{({\mathfrak g},{\mathfrak g}_{\0})}(M)$ and 
${\mathcal V}_{({\mathfrak e},{\mathfrak e}_{\0})}(M)$.  As \eqref{E:introiso} suggests, there is a close relationship between these two varieties.  More precisely, in Theorem~\ref{T:genericequality} we prove that ${\mathcal V}_{({\mathfrak g},{\mathfrak g}_{\0})}(M)$ is ``generically'' isomorphic to ${\mathcal V}_{({\mathfrak e},{\mathfrak e}_{\0})}(M)/W.$

Motivated by this relationship we more closely study the variety ${\mathcal V}_{({\mathfrak e},{\mathfrak e}_{\0})}(M)$ in Section~\ref{S:supportvarieties}.  We prove that 
${\mathcal V}_{({\mathfrak e},{\mathfrak e}_{\0})}(M)$ can be identified 
via a ``rank variety'' description as a certain affine subvariety of 
${\mathfrak e}_{\1}$. This concrete description allows us to verify that the operator ${\mathcal V}_{({\mathfrak e},{\mathfrak e}_{\0})}(-)$ satisfies many of the desirable properties of a support variety theory.  In particular, it detects projectivity in the sense of Dade's Lemma for finite groups and has the tensor product property (see Theorems~\ref{L:tensorproduct} and~\ref{C:tensorproducttheorem}).   Furthermore, one sees that the representation theory for the superalgebra ${\mathfrak e}$ over ${\mathbb C}$ 
has similar features to that of modular representations over 
fields of characteristic two (cf.\ Corollary~\ref{C:rankcor1}).  One interesting facet of our 
definition of rank varieties is that it involves restricting 
supermodules to copies of the Lie superalgebra ${\mathfrak q}(1)$ rather 
than cyclic shifted subgroups.

%In the process of our investigation
%we formulate several interesting conjectures relating the support varieties and relative 
%cohomology theories of $({\mathfrak g},{\mathfrak g}_{\0})$ and 
%$({\mathfrak e},{\mathfrak e}_{\0})$. These conjectures are interrelated with  
%earlier invariant theory results due to Panyushev \cite{panyushev}.  

In using these geometric and cohomological methods we also aim to obtain a  
deeper understanding of the combinatorics of the finite-dimensional 
representations of the Lie superalgebra ${\mathfrak g}$.  Let $\g$ be a simple classical Lie superalgebra with a nondegenerate invariant supersymmetric even bilinear form.  Kac and Wakimoto \cite{kacwakimoto} use the form to give a combinatorial definition of the (numerical) defect of $\g.$  The subsuperalgebras we introduce can be viewed as ``defect subalgebras'' of the principal block whose odd dimension coincides with the defect described above.  We also provide a cohomological interpretation of the defect of a Lie superalgebra which suggests a natural extension of the notion to Lie superalgebras without such a bilinear form.

Kac and Wakimoto also use the form on $\g$ to define the atypicality of a block and of a simple $\g$-supermodule. Atypicality is a combinatorial 
invariant used to give a rough measure of the complexity of these objects.  For instance, the characters of the typical (i.e.\ atypicality zero) simple supermodules for the basic classical Lie superalgebras were obtained early on by Kac \cite{Kacnote}.  On the other hand, much effort has gone into studying atypical simple supermodules, 
with most investigations being handled on a case by case basis (e.g.\  \cite{MVdJ, serganova2, serganova1, VdJZ, VdJ1}). 
In the last section of the paper we discuss how both the defect and 
atypicality are related to our support variety constructions.
We formulate an intriguing conjecture regarding this relationship which suggests a functorial description of atypicality.  We also discuss how our conjecture is closely related to a conjecture of Kac and Wakimoto.

%This provides ample evidence that these constructions will naturally allow one to 
%extend the notion of atypicality to all finite dimensional ${\mathfrak g}$-supermodules in a 
%functorial way. 

We remark that the approach taken here differs considerably 
from that of Duflo and Serganova \cite{dufloserganova}.  Certainly our varieties are not equal to theirs as our rank varieties contain elements whose bracket with themselves is non-zero.  As far as we know, the associated varieties they construct do not have a cohomological interpretation.  It would be of great interest to describe the relationship, if any, between these two theories.  

The authors would like to acknowledge David Benson, William Graham, Markus Hunziker, Gerald Schwarz, and Robert Varley 
for sharing their insights during various stages of this project.

\section{Lie Superalgebras, Representations, and Cohomology} \Label{S:cohom}

\subsection{Lie Superalgebras and their Representations}\Label{S:prelims}\setcounter{Df}{0}  Throughout 
we work with the complex numbers $\C$ as the ground field.  Recall that a superspace is a $\Z_{2}$-graded 
vector space and, given a superspace $V$ and a homogeneous vector $v \in V,$ we write $\p{v} \in \Z_{2}$ 
for the \emph{parity} (or \emph{degree}) of $v.$  Elements of $V_{\0}$ (resp. $V_{\1}$) are called \emph{even} 
(resp. \emph{odd}).  Note that if $M$ and $M'$ are two superspaces, then the space $\Hom_{\C}(M,M')$ is 
naturally $\Z_{2}$-graded by declaring $f \in \Hom_{\C}(M,M')_{r}$ ($r\in \Z_{2}$) if $f(M_{s}) \subseteq M'_{s+r}$ for 
all $s \in \Z_{2}$.   

A superalgebra is a $\Z_2$-graded, unital, associative algebra $A=A_{\0}\oplus A_{\1}$ which satisfies 
$A_{r}A_{s}\subseteq A_{r+s}$ for all $r,s\in \Z_2.$  A \emph{Lie superalgebra} is a superspace 
$\g=\g_{\0}\oplus \g_{\1}$ with a bracket operation $[\;,\;]:\g \otimes \g \to \g$ which preserves the 
$\Z_{2}$-grading and satisfies graded versions of the usual Lie bracket axioms.  In particular, we note 
that $\g_{\0}$ is a Lie algebra under the bracket obtained by restricting the bracket of $\g.$  If $\g$ is a 
Lie superalgebra, then one has a universal enveloping superalgebra $U(\g)$ which is $\Z_2$-graded and 
satisfies a PBW type theorem.  See, for example, \cite{K} for details and further background on Lie superalgebras.

We call a finite dimensional Lie superalgebra \emph{classical} if there is a connected reductive algebraic group $G_{\0}$ such that $\operatorname{Lie}(G_{\0})=\g_{\0},$ and an action of $G_{\0}$ on $\g_{\1}$ which differentiates to the adjoint action of $\g_{\0}$ on $\g_{\1}.$  In particular, if $\g$ is classical, then $\g_{\0}$ is a reductive Lie algebra and $\g_{\1}$ is semisimple as a $\g_{\0}$-module.  Note that we do not assume that $\g$ is simple.  A \emph{basic classical} Lie superalgebra is a classical Lie superalgebra with a nondegenerate invariant supersymmetric even bilinear form.  The simple (basic) classical Lie superalgebras were classified by Kac \cite{K}.

Given a Lie superalgebra, $\g$, let us describe the category of $\g$-supermodules.  The objects are 
all left $U(\g)$-modules which are $\Z_{2}$-graded; that is, superspaces $M$ satisfying $U(\g)_{r}M_{s} 
\subseteq M_{r+s}$ for all $r, s \in \Z_{2}.$  If $M$ is a $\g$-supermodule, then by definition $N \subseteq M$ is a 
subsupermodule if it is a supermodule which inherits its grading from $M$ in the sense that $M_{r} 
\cap N = N_{r}$ for $r \in \Z_2$. Given ${\mathfrak g}$-supermodules $M$ and $N$ one can use the 
antipode and coproduct of $U({\mathfrak g})$ to define a ${\mathfrak g}$-supermodule 
structure on the contragradient dual $M^{*}$ and the tensor product $M\otimes N$. 

A morphism of $U(\g)$-supermodules is an element of $\text{Hom}_{{\mathbb C}}(M,M^{\prime})$ 
satisfying $f(xm)=(-1)^{\p{f}\;\p{x}}xf(m)$ for all $m \in M$ and 
all $x \in U(\g).$  Note that this 
definition makes sense as stated only for homogeneous elements; it should be interpreted via linearity 
in the general case.  We emphasize that we allow \emph{all} morphisms and not just graded (i.e.\ \emph{even}) 
morphisms.  However, note that  $\text{Hom}_{U(\g)}(M,M^{\prime})$ inherits a $\Z_{2}$-grading as a subspace of  $\text{Hom}_{{\mathbb C}}(M,M^{\prime}).$ 

The category of $\g$-supermodules is not an abelian category.  However, the \emph{underlying even category}, 
consisting of the same objects but only the even morphisms, is an abelian category.  This, along with the parity 
change functor, $\Pi,$ which simply interchanges the $\Z_{2}$-grading of a supermodule, allows one to make use of 
the tools of homological algebra.

As a special case of the above discussion, we always view a Lie algebra (e.g.\ the even part of a Lie superalgebra) 
as a Lie superalgebra concentrated in degree $\0$.   

There are two categories of $\g$-supermodules which are natural to consider.  First, recall that one 
says a supermodule $M$ is {\em finitely semisimple} if it is isomorphic to a direct sum of finite dimensional 
simple subsupermodules.  If $\g$ is a Lie superalgebra and $\mathfrak{t}$ is a Lie subsuperalgebra, then 
let $\mathcal{C}=\mathcal{C}_{(\g,\mathfrak{t})}$ denote the full subcategory of the category of all $\g$-supermodules obtained by taking as objects all $\g$-supermodules which are finitely semisimple as $\mathfrak{t}$-supermodules.

We let $\mathcal{F}=\mathcal{F}_{(\g,\mathfrak{t})}$ denote the full subcategory of $\mathcal{C}$ obtained by taking the objects to be all finite dimensional supermodules in $\mathcal{C}.$  Observe that in the special case when $\mathfrak{t}=\g_{\0}$ is semisimple as a Lie algebra, then $\mathcal{F}$ is simply the category of finite dimensional $\g$-supermodules.  

As discussed in \cite[3.1.6]{kumar}, the category $\mathcal{C}$ is closed under arbitrary direct sums, quotients, and finite tensor products.

\subsection{Relative Cohomology}\Label{SS:relcohom} \setcounter{Df}{0} In this subsection we outline the basic 
definitions and results for relative cohomology for Lie superalgebras.  Relative cohomology for 
Lie algebras was first defined by Hochschild \cite{Hoch} and the super case is considered in Fuks \cite{fuks}.  
The main theme is that, once one accounts for the ${\mathbb Z}_2$-grading, results from the purely 
even case hold here as well.  For the sake of brevity we omit proofs when they are  
straightforward generalizations of the classical arguments.  We refer the reader to \cite[Appendix D]{kumar} 
for the details.

Let $R$ be a superalgebra and $S$ a subsuperalgebra.  In particular, we assume $S_{r}=R_{r} \cap S$ for $r \in \Z_2.$  Let 
\[
\dotsb \to M_{i-1} \xrightarrow{f_{i-1}} M_{i} \xrightarrow{f_{i}}  M_{i+1} \to \dotsb 
\]
be a sequence of $R$-supermodules and even $R$-supermodule 
homomorphisms.  We say this sequence is \emph{$(R,S)$-exact} if it is exact as a 
sequence of $R$-supermodules and if, when viewed as a sequence of $S$-supermodules,  
$\Ker f_{i}$ is a direct summand of $M_{i}$ for all $i.$  Note that our assumption 
that $f_{i}$ is even implies that $\Ker f_{i}$ is a subsupermodule of $M_{i}$ and that the 
splitting $M_{i}= \Ker f_{i} \oplus N_{i}$ is as $S$-supermodules.

An $R$-supermodule $P$ is \emph{$(R,S)$-projective} if given any $(R,S)$-exact sequence 
\[
0 \to M_{1} \xrightarrow{f} M_{2} \xrightarrow{g}  M_{3} \to 0,
\] 
and $R$-supermodule homomorphism $h: P \to M_{3}$ there is a $R$-supermodule map $\tilde{h}:P \to M_{2}$ satisfying 
$g \circ \tilde{h}=h.$  

In particular, if $P$ is a projective $R$-supermodule, then it is automatically $(R,S)$-projective.  Also, 
since $g$ is assumed to be even, if $h$ is homogeneous, then one can choose $\tilde{h}$ to be homogeneous of the same 
degree as $h.$  More generally, if we write $h=h_{\0}+h_{\1}$ where $h_{r} \in \Hom_{R}(P,M_{3})_{r}$ 
($r =\0 ,\1$), then we can lift each $h_{r}$ and $\widetilde{h_{\0}}+\widetilde{h_{\1}}$ is a lift of $h.$

An $(R,S)$-projective resolution of an $R$-supermodule $M$ is an $(R,S)$-exact sequence 
\[
\dotsb \xrightarrow{\delta_{2}} P_{1} \xrightarrow{\delta_{1}} P_{0} \xrightarrow{\delta_{0}} M \to 0,
\]
where each $P_{i}$ is an $(R,S)$-projective supermodule.  We remind the reader that implicit in the 
definition is the fact that the maps $\delta_{i}$ are all assumed to be even.

The following lemma is proven just as in \cite{kumar}.

\begin{lemma}\Label{L:projectives} Let $R$ be a superalgebra and $S$ be a subsuperalgebra of $R$.   
\begin{itemize} 
\item[(a)] If $M$ is any $S$-supermodule, then $R \otimes_{S} M$ is an  
$(R,S)$-projective $R$-supermodule.  The $\Z_{2}$-grading on $R \otimes_{S} M$ is given in the usual way by 
\[
(R \otimes_{S} M)_{i} =\bigoplus_{\substack{k,l \in \Z_{2}\\ k+l=i}} R_{k} \otimes_{S} M_{l}.
\]
\item[(b)] Any $R$-supermodule $M$ admits an $(R,S)$-projective resolution.  Namely,
\[
\dotsb \xrightarrow{\delta_{2}} R \otimes_{S} \Ker \delta_{0} \xrightarrow{\delta_{1}} R \otimes_{S}M \xrightarrow{\delta_{0}} M\to 0.
\] Here $\delta_{i}$ is the ``multiplication'' map $R \otimes_{S} N \to N$ given by $r \otimes n \mapsto rn$ for any $R$-supermodule $N.$ 
\end{itemize}
\end{lemma}

Note that since the multiplication map is even, its kernel is an $R$-subsupermodule of the domain and, hence, we can 
recursively define the above sequence as indicated. Given an $R$-supermodule $M$ with $(R,S)$-projective resolution 
$P_{\bullet} \to M,$ apply the functor $\Hom_{R}(-, N)$ and set 
\[
\Ext^{i}_{(R,S)}(M,N)=\operatorname{H}^{i}(\Hom_{R}(P_{\bullet}, N)).
\]  
One can show that $\text{Ext}^{\bullet}_{(R,S)}(M,N)$ can also be obtained using the dually defined 
$(R,S)$-injective resolutions, and that it 
is functorial in both arguments and well defined.  See \cite{kumar} for details. It is of interest to track 
the $\Z_2$-grading. As we remarked earlier, $\Hom_R(P_{i}, N)$ is naturally $\Z_2$-graded and since $\delta_{i}$ was 
assumed to be even the 
induced homomorphism $\Hom_R(P_{i-1}, N) \to \Hom_R(P_{i}, N)$ is also even.  Consequently, $\Ext^{i}_{(R,S)}(M,N)$ 
inherits a $\Z_2$-grading.

\subsection{Relative Cohomology for Lie Superalgebras}\Label{SS:lierelcohom} \setcounter{Df}{0}Let $\g$ be a Lie superalgebra and let $\mathfrak{t} \subseteq \g$ be a Lie subsuperalgebra.  
In this section we define (relative) Lie superalgebra cohomology for the pair $(\g,\mathfrak{t})$ and prove that it 
coincides with the relative cohomology defined in the previous section for the enveloping superalgebras of the pair. 

First let us recall the definition of Lie superalgebra cohomology.  Let $M$ be a $\g$-supermodule.  
For $p\geq 0$, set 
\[
C^{p}(\g,M)=\Hom_{{\mathbb C}}(\wedge^{p}_{s}(\g),M),
\] 
where $\wedge^{p}_{s}(\g)$ is the \emph{super} wedge product.  That is, $\wedge^{p}_{s}(\g)$ is the 
$p$-fold tensor product of $\g$ modulo the $\g$-subsupermodule generated by elements of the form 
\[
x_{1} \otimes \dotsb\otimes  x_{k} \otimes x_{k+1} \otimes \dotsb \otimes x_{p} + (-1)^{\p{x}_{k}\p{x}_{k+1}}
x_{1} \otimes \dotsb\otimes  x_{k+1} \otimes x_{k} \otimes \dotsb \otimes x_{p},
\] for homogeneous $x_{1}, \dotsc , x_{p} \in \g.$  Thus $x_{k},x_{k+1}$ skew commute unless both are odd in which case they commute.

Let 
\[
d^{p}: C^{p}(\g,M) \to C^{p+1}(\g,M)
\]
be given by 
\begin{multline}\label{E:differential}
d^{p}(\phi)(x_{1}\wedge \dotsb \wedge x_{p+1})\\
=\sum_{i < j} (-1)^{\sigma_{i,j}(x_{1}, \dotsc , x_{p})} 
\phi([x_{i},x_{j}] \wedge x_{1} \wedge \dotsb \wedge \hat{x}_{i}\wedge \dotsb \wedge \hat{x}_{j}\wedge \dotsb  \wedge x_{p+1}) \\
+ \sum_{i}(-1)^{\gamma_{i}(x_{1}, \dotsc , x_{p},\phi)} x_{i}\phi (x_{1} \wedge \dotsb \wedge \hat{x}_{i} 
\wedge \dotsb \wedge x_{p+1}),
\end{multline}
where $x_{1}, \dotsc , x_{p+1}$ and $\phi$ are assumed to be homogeneous, and   
\begin{align*}
\sigma_{i,j}(x_{1}, \dotsc , x_{p})&:=i+j+\p{x}_{i}(\p{x}_{1}+\dotsb +\p{x}_{i-1})+\p{x}_{j}(\p{x}_{1}+\dotsb +\p{x}_{j-1}+\p{x}_{i}),\\
\gamma_{i}(x_{1}, \dotsc , x_{p}, \phi)&:=i+1+\p{x}_{i}(\p{x}_{1}+\dotsb + \p{x}_{i-1}+\p{\phi}).
\end{align*}

Then we define 
\[
\operatorname{H}^{p}(\g, M)=\Ker d^{p}/\operatorname{Im} d^{p-1}.
\]

Now consider the relative version of the above construction.  Let $\g$, ${\mathfrak t}$, and $M$ be as 
above. Define 
\[
C^{p}(\g, \mathfrak{t}; M)=\Hom_{\mathfrak{t}}(\wedge^{p}_{s}(\g/\mathfrak{t}), M).
\]  Then the map $d^{p}$ descends to give a map $d^{p}: C^{p}(\g, \mathfrak{t}; M) \to  C^{p+1}(\g, \mathfrak{t}; M)$ and we define 
\[
\operatorname{H}^{p}(\g, \mathfrak{t}; M)=\Ker d^{p}/\operatorname{Im} d^{p-1}.
\]

As discussed earlier $\Hom$ spaces are naturally $\Z_{2}$-graded so the cochains have a $\Z_{2}$-grading.  
Note that the map $d^{p}$ preserves this grading and so the cohomology inherits a $\Z_{2}$-grading.  Thus 
cohomology is $\Z \times \Z_{2}$-graded --- by the cohomological degree and the $\Z_{2}$-grading, respectively.

Let us now relate the two cohomologies we have introduced.  The proof given in \cite{kumar} can be 
adapted to the super setting to yield the following result.  

\begin{lemma}\Label{L:relatescohom2} Let $\g$ be a Lie superalgebra, $\mathfrak{t}$ a Lie subsuperalgebra, and 
$M,N$ $\g$-supermodules.  Assume that $\g$ is finitely semisimple as a $\mathfrak{t}$-supermodule under 
the adjoint action. Then
\begin{equation*}
\Ext^{\bullet}_{(\mathcal{U}(\g),\mathcal{U}(\mathfrak{t}))}(M, N) \cong 
\Ext^{\bullet}_{(\mathcal{U}(\g),\mathcal{U}(\mathfrak{t}))}(\C, M^{*} \otimes N)
\cong \operatorname{H}^{\bullet}(\g, \mathfrak{t}; M^{*} \otimes N).
\end{equation*}
\end{lemma}

\subsection{Categorical Cohomology.} \Label{SS:cohomII}\setcounter{Df}{0} In this subsection we interpret 
relative cohomology in terms of cohomology in certain categories.  Our approach is inspired by that of \cite{gabber}.  
Fix a Lie superalgebra $\g$ and a Lie subsuperalgebra $\mathfrak{t}$.  Let 
$\mathcal{C}=\mathcal{C}_{(\g,\mathfrak{t})}$ and $\mathcal{F}=\mathcal{F}_{(\g,\mathfrak{t})}$ be the 
categories introduced in Section~\ref{S:prelims}. 

\begin{prop}\Label{T:finitelyss}  Assume $\g$ is finitely semisimple as a $\mathfrak{t}$-supermodule 
under the adjoint action.  
Let $L$ be a finite dimensional simple $\mathfrak{t}$-supermodule.  Let 
\[
\widehat{L} = \mathcal{U}(\g) \otimes_{\mathcal{U}(\mathfrak{t})} L.
\] 
Then $\widehat{L}$ is $(\mathcal{U}(\g),\mathcal{U}(\mathfrak{t}))$-projective and a projective supermodule 
in the category $\mathcal{C}_{(\g,\mathfrak{t})}.$
\end{prop}

\begin{proof} Note that in the definition of $\widehat{L}$ we view $\mathcal{U}(\g)$ as a right 
$\mathcal{U}(\mathfrak{t})$-supermodule via right multiplication, and $\widehat{L}$ as a 
left $\mathcal{U}(\g)$-supermodule 
via left multiplication.  By Lemma~\ref{L:projectives} $\widehat{L}$ is 
$(\mathcal{U}(\g),\mathcal{U}(\mathfrak{t}))$-projective.

Now since $\mathcal{C}$ is closed under tensor products, arbitrary direct sums, and quotients, $\mathcal{U}(\g)$ is 
a finitely semisimple $\mathfrak{t}$-supermodule under the \emph{adjoint} action.  
Consequently, $\mathcal{U}(\g) \otimes_{\C} L$ is a finitely semisimple $\mathfrak{t}$-supermodule. 
It remains to observe that the induced module $\widehat{L}$ can be obtained as a 
quotient of $\mathcal{U}(\g) \otimes_{\C} L$. (cf.\ the proof of \cite[Corollary 3.1.8]{kumar}).  
Consequently, $\widehat{L}$ is finitely semisimple, hence is an object in $\mathcal{C}_{(\g,\mathfrak{t})}.$
Finally, one can verify that $\widehat{L}$ is a projective object in $\mathcal{C}_{(\g,\mathfrak{t})}$ 
using Frobenius reciprocity. 
\end{proof}

\begin{corollary}\Label{C:extequivalence}Assume $\g$ is finitely semisimple as a $\mathfrak{t}$-supermodule under 
the adjoint action.  Let $\mathcal{D}$ be a full subcategory of $\mathcal{C}_{(\g,\mathfrak{t})}$ such that whenever $L$ 
is a finite dimensional simple $\mathfrak{t}$-supermodule which appears as a composition factor of some object in 
$\mathcal{D},$ then $\widehat{L}$ is an object in $\mathcal{D}.$  Furthermore, assume $\mathcal{D}$ is closed under direct sums.
Then $\mathcal{D}$ has enough projectives and for any $M,N$ which are objects in $\mathcal{D}$ we have
\[
\Ext^{p}_{\mathcal{D}}(M,N) \cong \Ext^{p}_{(\mathcal{U}(\g),\mathcal{U}(\mathfrak{t}))}(M,N),
\] as superspaces.
\end{corollary}

\begin{proof}  Let $M$ be an object in $\mathcal{D}$ and write $M$ as the direct sum of finite dimensional 
simple $\mathfrak{t}$-supermodules:
\[
M= \bigoplus_{j \in I} L_{j}. 
\] Then by Frobenius reciprocity one sees that $P:=\bigoplus_{j \in I} \widehat{L}_{j}$ maps onto $M.$  
By the previous proposition $P$ is projective in $\mathcal{D}$ and so $\mathcal{D}$ has enough projectives.  
Therefore, we also see that any object in $\mathcal{D}$ has a resolution by supermodules which are both projective 
in $\mathcal{D}$ and $(\mathcal{U}(\g), \mathcal{U}(\mathfrak{t}))$-projective. In order to compute 
cohomology in both categories one applies $\text{Hom}_{U({\mathfrak g})}(-,N)$ to this resolution, 
thus the equivalence of $\Ext$s follows.
\end{proof}

\subsection{Cohomology for the pair $(\g,\g_{\0})$}\Label{SS:ggzerocohom}\setcounter{Df}{0}
\emph{For the remainder of the article we will make the following assumptions.  We will always assume $\g$ is a classical Lie superalgebra.  A $\g$-supermodule will always be assumed to be an object in the category $\mathcal{C}=\mathcal{C}_{(\g,\g_{\0})}$ and a finite dimensional $\g$-supermodule will always mean an object in the category $\mathcal{F}=\mathcal{F}_{(\g,\g_{\0})}.$} 

By the PBW theorem for Lie superalgebras (cf.\ \cite[1.1.3]{K}), if $N$ is a finite 
dimensional $\g_{\0}$-supermodule, then 
\[
U(\g) \otimes_{U(\g_{\0})} N
\] is again finite dimensional.  Therefore the conditions of Corollary~\CORextequivalence\ are satisfied 
for the pair $(\g,\g_{\0})$ where the category $\mathcal{D}$ is $\mathcal{F}.$   Combining that result with 
Lemma~\Lrelatescohom2\ yields the following theorem.

\begin{theorem} \Label{T:cohomforfdmoduls} Let $M$ and $N$ be finite dimensional $\g$-supermodules.  Then for all $p \geq 0$ we have 
\[
\Ext^{p}_{\mathcal{F}}(M,N) \cong \operatorname{H}^{p}(\g, \g_{\0}; M^{*} \otimes N).
\] 
\end{theorem}

In particular, we can apply the above result when $M=N=\C$ to obtain the cohomology ring for $\g$ in the category $\mathcal{F}.$  
Before doing so, let us make several observations which will simplify the calculation.

 We first note that the differential $d^{p}$ in 
\eqref{E:differential} is identically zero.  Namely, since the bracket preserves the $\Z_2$-grading we have $[\g_{\1},\g_{\1}] \subseteq \g_{\0}.$  From this observation one concludes that in 
the first sum of \eqref{E:differential} each $[x_{i}, x_{j}]$ is always zero in the quotient $\g/\g_{\0}$ and hence 
these terms are identically zero.  The terms in the second sum of \eqref{E:differential} are all zero since here $M$ is 
the trivial supermodule.  Consequently, the cohomology is simply the cochains themselves.  Second, we observe that $\g/\g_{\0} \cong \g_{\1}$ as a $\g_{\0}$-supermodule.  Taken together, these observations imply
\begin{align*}
\Ext^{p}_{\mathcal{F}}(\C,\C) &\cong \Hom_{\g_{\0}}(\wedge^{p}_{s}(\g_{\1}),\C).
\end{align*}

However, one can simplify further with the following observation.  Since $\g_{\1}$ is purely odd the elements 
of $\wedge^{p}_{s}(\g_{\1})$ commute without sign and so this super wedge product can be viewed as 
a classical symmetric product.  Let $G_{\0}$ denote the connected reductive group with Lie algebra $\g_{\0}.$  
We then obtain the following theorem.

\begin{theorem}\Label{T:cohomring} Let $\g$ be a classical Lie superalgebra and
let $\mathcal{F}$ be the category of finite dimensional $\g$-supermodules.  Then,
\[
\Ext^{\bullet}_{\mathcal{F}}(\C,\C) \cong \operatorname{H}^{\bullet}(\g,\g_{\0}; \C)  \cong S(\g_{\1}^{*})^{\g_{\0}} =S(\g_{\1}^{*})^{G_{\0}}.
\] 
\end{theorem}  The last equality is a well known result in characteristic zero Lie theory (cf.\ \cite[24.3.3]{tauvel}).  
As a matter of notation we write $S(\g_{\1}^{*})$ for $S^{\bullet}(\g_{\1}^{*}).$

Note that one can check by the definition that this isomorphism is one of superalgebras; that is, the map respects multiplication.  
Since $S(\g_{\1}^{*})$ is generated by elements of degree $(1, \1 )$ in the $\mathbb{Z} \times 
\mathbb{Z}_{2}$-grading, both it and the cohomology ring will consist of elements whose $\Z_{2}$ degree is the reduction modulo $2$ of the cohomological degree.  In particular, the cohomology ring is always a commutative ring in the ungraded sense.  Because of this and as the $\Z_{2}$-grading will not generally play an important role in what follows, we will will leave it implicit unless otherwise necessary.

Recall that if $M$ is a finite dimensional $\g$-supermodule, then  
$\Ext^{\bullet}_{\mathcal{F}}(\C, M)=\operatorname{H}^{\bullet}(\g,\g_{\0}; M)$ is a  
$\Ext^{\bullet}_{\mathcal{F}}(\C,\C)=\operatorname{H}^{\bullet}(\g,\g_{\0};\C)$-module via the Yoneda product 
(cf.\ \cite[5.7]{benson}).  

\begin{theorem} \Label{T:finitegen2}Let $\g$ be a classical Lie superalgebra.  Let $M$ be a finite dimensional $\g$-supermodule.  The superalgebra 
$\operatorname{H}^{\bullet}(\g,\g_{\0}; \C)$  is finitely generated as a ring.  
Furthermore,  $\operatorname{H}^{\bullet}(\g,\g_{\0};M)$ 
is finitely generated as an $\operatorname{H}^{\bullet}(\g,\g_{\0}; \C)$-module.
\end{theorem}

\begin{proof}  The fact that $\operatorname{H}^{\bullet}(\g,\g_{\0}; \C)$ is finitely generated is immediate from Theorem~\THcohomring\ and the 
classic invariant theory result of Hilbert \cite[Theorem 3.6]{pv}.  Now consider the second statement of the theorem.  Clearly, 
\[
\bigoplus_{p \geq 0} \Hom_{\C}(\wedge^{p}_{s}(\g_{\1}), M) \cong \bigoplus_{p \geq 0} \left( \wedge^{p}_{s}(\g_{\1}^{*}) 
\otimes M \right)\cong S(\g_{\1}^{*}) \otimes M
\] is finitely generated as an $S(\g_{\1}^{*})$-module.   Applying \cite[3.25]{pv}, we see that 
\[
\left( \bigoplus_{p \geq 0} \Hom_{\C}(\wedge^{p}_{s}(\g_{\1}), M) \right)^{\g_{\0}}=\bigoplus_{p \geq 0} 
\Hom_{\g_{\0}}(\wedge^{p}_{s}(\g_{\1}), M)
\] is finitely generated as a $S(\g_{\1}^{*})^{G_{\0}}=\operatorname{H}^{\bullet}(\g,\g_{\0}; \C)$-module.  
That is, the cochain complex $C^{\bullet}(\g, \g_{\0}; M)$ used in Section~\lierelcohom\ to define $\operatorname{H}^{\bullet}(\g,\g_{\0};M)$ is finitely generated as a module over $\operatorname{H}^{\bullet}(\g,\g_{\0}; \C).$  

The finite generation of $C^{\bullet}(\g, \g_{\0}; M)$ implies that $\operatorname{H}^{\bullet}(\g, \g_{\0}; M)$ is finitely generated as follows.  Given $r \in \operatorname{H}^{p}(\g,\g_{\0}; \C)$ and $x \in C^{\bullet}(\g,\g_{\0};M),$ one has $d(rx)=d(r)x+(-1)^{p}rd(x)=(-1)^{p}rd(x),$ where $p$ is the cohomological degree of $r$ (cf.\ the proof of \cite[Theorem 5]{Gru:97}).  The second equality follows from the fact that the differentials for  $\operatorname{H}^{\bullet}(\g, \g_{\0}; \C)$ are identically zero.  Therefore, $d:C^{\bullet}(\g, \g_{\0}; M) \to C^{\bullet}(\g, \g_{\0}; M)$ is a graded  $\operatorname{H}^{\bullet}(\g,\g_{\0}; \C)$-module homomorphism.  Since $\operatorname{H}^{\bullet}(\g,\g_{\0}; \C)$ is finitely generated, any subquotient of a finitely generated graded module is finitely generated.  This implies the final statement of the theorem.
\end{proof}

\section{Invariant Theory} \Label{S:alggeo}
%\vskip .25cm 
\subsection{}  \setcounter{Df}{0} In this section we develop some of the invariant theory needed to
study cohomology in the category $\mathcal{F}.$  
More specifically, we will study the action of $G_{\0}$ on $\g_{\1}$. When the 
action is either stable (see Section~\ref{SS:stableactions}) or polar (see Section~\ref{SS:dadokkac}) 
there are general results about rings of invariants which will be central to the approach taken here.  These results will be used in Section~\ref{S:morecohom} 
to provide more structural information pertaining to the relative cohomology groups for these 
Lie superalgebras.   
By doing a 
case by case analysis one can determine which simple classical Lie superalgebras and related algebras 
admit stable and/or polar actions.  This data is presented in Table~5 of the Appendix. 

We view $\g_{\1}$ as an affine 
variety with the Zariski topology.  It is isomorphic to the affine space $\mathbb{A}^{\operatorname{dim}\g_{\1}}$ and so is smooth and 
irreducible. The group $G_{\0}$ acts via the adjoint action.  Given $g \in G_{\0}$ and $x \in \g_{\1},$ we write $g.x$ 
for this action.  Given $x \in \g_{\1},$ let $G_{\0,x}$ denote the stabilizer subgroup 
$G_{\0,x}=\setof{g \in G_{\0}}{g.x=x}$ and let $G_{\0}.x=\{g.x \:\vert\:  g \in G_{\0} \}$.

\subsection{Stable Actions}\Label{SS:stableactions} \setcounter{Df}{0} Let us recall the following definitions and results from invariant theory.

\vskip .15cm 
\begin{itemize}
\item[(a)] A point $x \in \g_{\1}$ is \emph{regular} if $G_{\0}.x$ has the maximum possible dimension.  Equivalently, $x$ is regular if $G_{\0 ,x}$ has minimal dimension.  By \cite[Sec. 1.4]{pv} there is a dense open subset of $\g_{\1}$ consisting of regular points.
\item[(b)] A point $x \in \g_{\1}$ is \emph{semisimple} if the orbit $G_{\0}.x$ is closed in $\g_{\1}.$  
\item [(c)]  The action of $G_{\0}$ on $\g_{\1}$ is called \emph{stable} if there is an open dense subset of 
$\g_{\1}$ consisting of semisimple points.  By a theorem of Popov \cite{popov} the action of $G_{\0}$ on $\g_{\1}$ 
will be stable if and only if $\g_{\1}$ has regular semisimple points.
\item[(d)]  If there is an open subset of $\g_{\1}$ such 
that the stabilizer subgroups of any two points in this set are conjugate subgroups of $G_{\0}$, then 
the stabilizer of such a point is called a \emph{stabilizer in general position}.  By a 
theorem of Richardson \cite[Theorem 7.2]{pv}, such an open set exists in $\g_{\1}.$  
\item [(e)] If the action of $G_{\0}$ on $\g_{\1}$ is stable, then one has that 
$\g_{\1 }$ contains an open dense set of regular semisimple elements whose stabilizers are in general position.  
We call such points \emph{generic}.
\end{itemize}

For short we say \emph{$\g$ is stable} if the action of $G_{\0}$ on $\g_{\1}$ is stable.  Assume $\g$ is stable.  Fix a generic element $x_{0} \in \g_{\1}$ and set
\[
H=G_{\0 , x_{0}}.
\]
Note that if $\pi: \g_{\1} \to \g_{\1}/G_{\0}$ is the canonical projection, then $\pi(x_{0}) \in \g_{\1}/G_{\0}$ is 
principal in the sense of \cite{luna}.  Therefore we can utilize the generalization of the 
Chevalley Restriction Theorem given by Luna and Richardson in \cite{luna}.  Set 
\[
\e_{\1}= \g_{\1}^{H}=\setof{z \in \g_{\1}}{h.z=z \text{ for all $h \in H$} },
\] and set 
\[
N=N_{G_{\0}}(H)=\setof{g \in G_{\0}}{gHg^{-1} = H}=\setof{g \in G_{\0}}{g.\e_{\1} = \e_{\1 }}.
\]  Since $x_{0}$ is semisimple one knows that $H$ is reductive and hence, by \cite[Lemma 1.1]{luna}, so is $N.$

We are now prepared to state the following key theorem.  

\begin{theorem} \Label{T:lunarichardson} Let $\g$ be a classical Lie superalgebra which is stable.  Let $M$ be a finite dimensional $G_{\0}$-module.
\begin{itemize}
\item[(a)] The restriction homomorphism $S(\g_{\1}^{*}) \to S(\e_{\1}^{*})$ induces an isomorphism 
\[
\operatorname{res}: S(\g_{\1}^{*})^{G_{\0}} \to S(\e_{\1}^{*})^{N}.
\]
\item[(b)] The set $G_{\0}.\e_{\1}$ is dense in $\g_{\1}.$ 
\item[(c)] The map induced by restriction 
\[
\rho: \Hom_{G_{\0}}(S^{n}(\g_{\1}),M) \to \Hom_{N}(S^{n}(\e_{\1}),M)
\] is injective for any $n \geq 0$.
\end{itemize} 
\end{theorem}

\begin{proof} (a--b). Part (a) is a direct application of \cite[Corollary 4.4]{luna}, and 
(b) follows by \cite[Section 7.1]{pv}.

(c). It is convenient to dualize the situation.  That is, we instead show \[
\rho: \Hom_{G_{\0}}(M, S^{n}(\g_{\1}^{*})) \to \Hom_{N}(M, S^{n}(\e_{\1}^{*}))
\] is injective for any $G_{\0}$-module $M$.  This, of course, is equivalent to the claimed statement.

Let $\varphi \in \Hom_{G_{\0}}(M, S^{n}(\g_{\1}^{*}))$ such that $\rho(\varphi)=0.$  That is, 
for any $m \in M,$ we have $\rho(\varphi)(m) \in S^{n}(\e_{\1}^{*})$ is zero as a homogenous 
degree $n$ polynomial on $\e_{\1}.$  In other words, $\varphi(m)(y)=0$ for all $m \in M$ and all 
$y \in \e_{\1}.$  Since $\varphi$ is a $G_{\0}$-homomorphism, for any fixed $m \in M$ we have 
\[
\varphi(m)(g.y)=\varphi(g^{-1}.m)(y)=0
\] for all $g \in G_{\0}$ and $y \in \e_{\1}.$  That is, $\varphi(m)$ is identically zero on $G_{\0}.\e_{\1}.$  
However, by part (b) this implies that $\varphi(m)=0.$  Since $m \in M$ was arbitrary, this implies 
that $\varphi$ itself is identically zero.
\end{proof}

\subsection{Polar Representations}\Label{SS:dadokkac} \setcounter{Df}{0} We now recall the notion of \emph{polar representations} introduced by 
Dadok and Kac \cite{dadokkac}.  Let $G$ be a reductive algebraic group acting on a vector space $V.$  Let $v \in V$ 
be a semisimple element.  Let 
\begin{equation}\Label{E:cv}
\c_{v}=\left\{x \in V \:\vert\:  \g.x \subseteq \g.v \right\},
\end{equation}
where $\g$ is the Lie algebra of $G.$  In general one has $\operatorname{dim} \c_{v} \leq \operatorname{dim} S(V^{*})^{G}.$  
By definition the action of $G$ on $V$ is \emph{polar} if for some semisimple $v \in V$ we have 
$\operatorname{dim} \c_{v} = \operatorname{dim} S(V^{*})^{G}.$  In this case $\c_{v}$ is called a \emph{Cartan subspace}.  For brevity, we say $\g$ is polar when the action of $G_{\0}$ on $\g_{\1}$ is polar.  We write $\c_{\1}$ for our fixed choice of a Cartan subspace.

In the case when the action of $G_{\0}$ on $\g_{\1}$ is both stable and polar, then one can further assume  
\begin{equation}\label{E:containments}
x_{0} \in \c_{x_{0}}=\c_{\1} \subseteq \e_{\1},
\end{equation}
where $x_{0}$ and $\e_{\1}$ are as in Section~\SSLR\ (cf.\ \cite[p. 514]{dadokkac}).  When the action is stable and polar, the Cartan subspace is unique up to conjugation by $G_{\0}$ by \cite[Theorem 2.3]{dadokkac}.

Let us set the following notation.  Given a subspace $V \subseteq \g_{\1},$ let 
\begin{align*}
\operatorname{Norm}_{G_{\0}}(V) &= \setof{g \in G_{\0}}{g.V = V},\\
\operatorname{Stab}_{G_{\0}}(V) &= \setof{g \in G_{\0}}{g.v = v\  \text{ for all $v \in V$}}. 
\end{align*} One has the following theorem.

\begin{theorem}\Label{T:dadokkac} Assume $\g$ is a classical Lie superalgebra which is polar.  Let $\c_{\1}$ be a Cartan subspace of $\g_{\1}.$   Let $M$ be a $G_{\0}$-module. Then
\begin{itemize}
\item[(a)] Restriction of functions then defines an isomorphism
\[
S(\g_{\1}^{*})^{G_{\0}} \cong S(\c_{\1}^{*})^{W},
\] where 
\[
W:=\operatorname{Norm}_{G_{\0}}(\c_{\1})/\operatorname{Stab}_{G_{\0}}(\c_{\1}).
\]  Furthermore, $W$ is necessarily a finite pseudoreflection group.
\item[(b)] If $\g$ is also stable, then $G_{\0}.\c_{\1}$ is dense in $\g_{\1}.$ 
\item[(c)] If $\g$ is also stable, then the map induced by restriction 
\[
\rho: \Hom_{G_{\0}}(S^{n}(\g_{\1}),M) \to \Hom_{N'}(S^{n}(\c_{\1}),M)
\] is injective for any $n \geq 0$, where $N':=\operatorname{Norm}_{G_{\0}}(\c_{\1})$. 
\end{itemize} 
\end{theorem}

\begin{proof} (a) This follows from \cite[Theorems 2.9--2.10]{dadokkac}. (b) Since we assume the action of $G_{\0}$ on $\g_{\1}$ is stable 
there is an open dense subset of $\g_{\1}$ consisting of semisimple elements.  By \cite[Proposition 2.2]{dadokkac} 
every semisimple element of $\g_{\1}$ is conjugate to some element of $\c_{\1}.$  The result then follows. (c) 
This follows by the same reasoning as in the proof of Theorem~\THlunarichardson(c). 
\end{proof}

Assume $\g$ is both stable and polar.  Recall that the \emph{discriminant} of the action of $W$ on $\c_{\1}$ is an element $D \in S(\c_{\1}^{*})^{W}$  where, for any $y \in \c_{\1},$  $D(y) \neq 0$ if and only if $W_{y}$ is trivial.  To be concrete, one can take $D$ to be a sufficiently high power of the Jacobian, $J$ (recalling that $J\in S(\c_{\1}^{*})$ is a skew invariant which satisfies $J(y) \neq 0$ if and only if $W_{y}$ is trivial).  The Jacobian for the polar Lie superalgebras considered in Sections~\ref{SS:typeA}--\ref{SS:typeF4} can be found in Table~5.  By doing a case by case check one can verify that for the simple classical Lie superalgebras which are stable and polar, $y \in \g_{\1}$ is generic if and only if $D(y) \neq 0.$  Therefore, in these cases the set of generic points of $\g_{\1}$ is precisely the principal open set defined by $D.$

\subsection{}\Label{SS:kostantrallis}\setcounter{Df}{0} We now recall a generalization of a theorem of Kostant and Rallis as presented in \cite[12.4.6]{goodman}. 

\begin{theorem} \label{T:KR} Let $K$ be a connected, reductive, linear algebraic group and let $V$ be a regular representation of $K$ and assume 
there is a subspace $\mathfrak{a} \subseteq V$ such that the following holds:
\begin{enumerate}
\item Restriction of functions defines an isomorphism of $S(V^{*})^{K}$ onto a subalgebra $\mathcal{R}$ of $S(\mathfrak{a}^{*}).$
\item The subalgebra $\mathcal{R}$ is generated by algebraically independent homogeneous elements $u_{1},\dotsc ,u_{r}$ with 
$r=\operatorname{dim} \mathfrak{a}.$  Furthermore, there exists a graded subspace $\mathcal{A}$ of $S(\mathfrak{a}^{*})$ such that 
the multiplication map $\mathcal{A} \otimes \mathcal{R} \to S(\mathfrak{a}^{*})$ is a linear isomorphism.
\item There exists $h \in \mathfrak{a}$ such that $|K.h \cap \mathfrak{a}| \geq \operatorname{dim} \mathcal{A}.$
\item Let $h$ be as above, and set 
\[
\mathcal{X}_{h}=\left\{v \in  V \:\vert\:  f(v)=f(h) \text{ for all } f \in S(V^{*})^{K}\right\}.
\]  If $v \in \mathcal{X}_{h},$ then 
\begin{equation}\Label{E:dimequality}
\operatorname{dim} K.v = \operatorname{dim} V - \operatorname{dim} \mathfrak{a}.
\end{equation}
\end{enumerate}  Assuming the above hypotheses hold, then there exists a graded $K$-submodule $\mathcal{H} \subseteq S(V^{*})$ such that multiplication provides an isomorphism of graded right $S(V^{*})^{K}$-modules, 
\[
\mathcal{H} \otimes S(V^{*})^{K} \to S(V^{*}).
\]  Furthermore, 
\[
\mathcal{H} \cong \Ind_{K_{h}}^{K} \C
\] as $K$-modules, where $K_{h}$ is the stabilizer of $h$ in $K.$

\end{theorem}

\subsection{} \Label{SS:KRforg}\setcounter{Df}{0} Throughout this subsection we assume $\g$ is stable and polar.  We first apply Theorem~\THKosRal\ to the action of $G_{\0}$ on $\g_{\1}.$  Set $K=G_{\0}$ and $V=\g_{\1}.$  Let $\mathfrak{a}=\c_{\1},$ the Cartan subspace fixed in Section~\ref{SS:dadokkac}.

\begin{theorem} \Label{T:KRforg} Let $x_{0} \in \g_{\1}$ be our fixed generic point and let $H=G_{\0, x_{0}}$.  Then there exists a graded $G_{\0}$-module $\mathcal{H} \subseteq S(\g_{\1}^{*})$ such that the multiplication map 
\[
\mathcal{H} \otimes S(\g_{\1}^{*})^{G_{\0}} \to S(\g_{\1}^{*})
\] gives an isomorphism of graded right $S(\g_{\1}^{*})^{G_{\0}}$-modules.  Furthermore, $\mathcal{H} \cong \Ind_{H}^{G_{\0}} \C$ as a $G_{\0}$-module.
\end{theorem}

\begin{proof}  We simply need to verify that the conditions of Theorem~\THKosRal\ apply.

\begin{enumerate}
\item By Theorem~\ref{T:dadokkac} restriction provides an isomorphism 
\[
 S(\g_{\1}^{*})^{G_{\0}} \to S(\c_{\1}^{*})^{W},
\] where $W$ is as given in the theorem.
\item Since $W$ is a finite pseudoreflection group, it follows immediately that $S(\g_{\1}^{*})^{G_{\0}}$ is generated by $\operatorname{dim} 
\c_{\1}$ algebraically independent generators.  There exists a graded subspace $\mathcal{A} \subseteq 
S(\c_{\1}^{*})$ such that multiplication $\mathcal{A} \otimes S(\c_{\1}^{*})^{W} \to S(\c_{\1}^{*})$ provides a graded 
linear isomorphism.  Furthermore, we have $\dim \mathcal{A}=|W|.$  See, for example, \cite[Theorem 12.4.7]{goodman} for a summary of these well-known results.  
\item Let $h=x_{0} \in \c_{\1}$. Since $\c_{\1} \subseteq \e_{\1}=\g_{\1}^{H},$ we have 
that $H \subseteq \operatorname{Stab}_{G_{\0}}(\c_{\1})$ and so we have 
\[
|G_{\0}.x_{0} \cap \c_{\1}| \geq |\operatorname{Norm}_{G_{\0}}(\c_{\1}).x_{0}| \geq |W.x_{0}|=|W|
=\operatorname{dim} \mathcal{A}.
\] 
\item If $\pi: \g_{\1} \to \g_{\1}/G_{\0}$ is the canonical quotient map, then by definition (see \cite[Sec. 4.4]{pv}) 
\[
\mathcal{X}_{x_{0}}=\pi^{-1}(\pi(x_{0})).
\]  However, since $x_{0}$ is semisimple and regular, the fiber $\pi^{-1}(\pi(x_{0}))$ is precisely $G_{\0}.x_{0}.$  
Therefore it suffices to check \eqref{E:dimequality} for $x_{0}.$  Now since $x_{0}$ is generic, by \cite[Theorem 7(ii)]{shafarevich} 
one has that 
\begin{align*}
\operatorname{dim} G_{\0}.x_{0}&=\operatorname{dim} \g_{\1}-\operatorname{dim} \g_{\1}/G_{\0}\\
             &=\operatorname{dim} \g_{\1}-\operatorname{dim} \c_{\1},
\end{align*} where the last equality follows from the fact that the Krull dimension of $S(\g_{\1}^{*})^{G_{\0}}$ equals the dimension of $\c_{\1}.$
\end{enumerate}
\end{proof}

Now consider the group $N$ acting on $\e_{\1}.$  Note that in the following theorem there is the additional assumption that the stabilizer of a generic point is connected.  However, we verified by direct calculation that this condition is satisfied for all the Lie superalgebras considered in Sections~\ref{SS:typeA}--\ref{SS:typeF4} which are stable and polar.  The stabilizer $H$ in these cases can be found in Table~4.

\begin{theorem} \label{T:KRfore} Let $x_{0} \in \g_{\1}$ be our fixed generic point and let $H=G_{\0, x_{0}}$.  Assume $H$ is connected. Then there exists a graded $N$-module $\mathcal{H}\subseteq S(\e_{\1}^{*})$ such that the multiplication map 
\[
\mathcal{H} \otimes S(\e_{\1}^{*})^{N} \to S(\e_{\1}^{*})
\] gives an isomorphism of graded right $S(\e_{\1}^{*})^{N}$-modules.  Furthermore, $\mathcal{H} \cong \Ind_{H}^{N} \C$ as $N$-modules.
\end{theorem}

\begin{proof}  Since $N$ is in general not connected one cannot apply Theorem~\ref{T:KR} directly.  Write $N^{0}$ for the connected component of the identity of $N$ and instead apply Theorem~\ref{T:KR} to $N^{0}$ acting on $\e_{\1}.$  This action is stable.  Indeed, one sees that $x_{0} \in \e_{\1}$ is a generic element as follows.  Since the $N^{0}$-stabilizer of any point of $\e_{\1}$ contains $H$ $x_{0}$ is regular.  To show that $x_{0}$ is semisimple, first note by Luna's criterion \cite[Theorem 6.17]{pv} that $N.x_{0}$ is closed in $\e_{\1}.$   Using that $H \subseteq N^{0},$ one sees that $N.x_{0}$ is the disjoint union of a finite number of $N^{0}$-orbits, each of the same dimension as $N^{0}.x_{0}.$  Therefore each of them, including $N^{0}.x_{0},$ is closed.

 The action of $N^{0}$ on $\e_{\1}$ is also polar with Cartan subspace $\c_{\1}$.  This can be seen as follows.  Since the action of $G_{\0}$ on $\g_{\1}$ is both stable and polar, \cite[Corollary 2.5]{dadokkac} implies $\g_{\1}=\c_{\1} \oplus \g_{\0}.x_{0}.$  Also, by assumption, $\c_{\1} \subseteq \e_{\1}.$  Taken together this implies $\c_{\1} \oplus \operatorname{Lie}(N^{0}).x_{0} \subseteq \e_{\1}.$  However, one has
\[
\dim \c_{\1}=\dim \g_{\1}/G_{\0 }=\dim \e_{\1}/N = \dim \e_{\1}/N^{0} = \dim \e_{\1}-\dim \operatorname{Lie}(N^{0}).x_{0}.
\] The first equality holds because $\c_{\1}$ is a Cartan subspace for the polar action of $G_{\0}$ on $\g_{\1},$  the second equality follows from Theorem~\ref{T:lunarichardson}(a), the third equality holds because $N/N^{0}$ is a finite group, and the last equality holds because $x_{0}$ is a generic point for the action of $N^{0}$ on $\e_{\1}$.
Therefore, by dimension counting, $\c_{\1} \oplus \operatorname{Lie}(N^{0}).x_{0} = \e_{\1}.$  This along with the fact that $\c_{\1}$ is a Cartan subspace for the action of $G_{\0}$ on $\g_{\1}$ implies that $\c_{\1}$ is a Cartan subspace for the action of $N^{0}$ on $\e_{\1}.$

Applying Theorem~\ref{T:KR} just as in the proof of Theorem~\ref{T:KRforg} one has that
\[
(\Ind_{H}^{N^{0}} \C) \otimes S(\e_{\1}^{*})^{N^{0}} \cong S(\e_{\1}^{*})
\] as graded right $S(\e_{\1}^{*})^{N^{0}}$-modules.  Since the action of $N^{0}$ on $\e_{\1}$ is polar, $S(\e_{\1}^{*})^{N^{0}}$ is a polynomial ring with an induced $N/N^{0}$ action.  Since $N/N^{0}$ is a finite group and 
\[
 (S(\e_{\1}^{*})^{N^{0}})^{N/N^{0}} \cong S(\e_{\1}^{*})^{N}
\] is a polynomial ring, it follows that $N/N^{0}$ is a pseudoreflection group.  By the classical theory of such groups \cite[Theorem 12.4.7]{goodman}, multiplication gives an isomorphism,
\[
\left(\Ind_{1}^{N/N^{0}}\C \right) \otimes \left(S(\e_{\1}^{*})^{N^{0}}\right)^{N/N^{0}} \cong S(\e_{\1}^{*})^{N^{0}},
\]  as graded right $S(\e_{\1}^{*})^{N}$-modules.  Therefore, combining these results, one has
\begin{align*}
S(\e_{\1}^{*}) &\cong \left(\Ind_{H}^{N^{0}} \C\right) \otimes S(\e_{\1}^{*})^{N^{0}} \\
                         &\cong \left(\Ind_{H}^{N^{0}} \C\right) \otimes \left(\Ind_{1}^{N/N^{0}}\C\right) \otimes S(\e_{\1}^{*})^{N},\\
                         &\cong \left( \Ind_{H}^{N^{0}} \C\right) \otimes \left(\Ind_{N^{0}}^{N}\C\right) \otimes S(\e_{\1}^{*})^{N}.
\end{align*}
The action of $N$ on  $\Ind_{1}^{N/N^{0}}\C$ is by inflation through the canonical map $N \to N/N^{0}.$  Under this action $\Ind_{1}^{N/N^{0}}\C \cong \Ind_{N^{0}}^{N}\C$ as $N$-modules, yielding the last isomorphism.  The action of $N \cong N^{0}\rtimes N/N^{0}$ on $\Ind_{H}^{N^{0}} \C$ is as described in \cite[I.3.8]{jantzen}.

By the tensor identity and transitivity of induction \cite[I.3.5-6]{jantzen} one obtains
\begin{align*}
S(\e_{\1}^{*}) &\cong \Ind_{N^{0}}^{N}\left(\Ind_{H}^{N^{0}} \C \otimes \C\right) \otimes S(\e_{\1}^{*})^{N},\\
     &\cong \left( \Ind_{N^{0}}^{N}\Ind_{H}^{N^{0}} \C\right) \otimes S(\e_{\1}^{*})^{N}\\
     & \cong \left(\Ind_{H}^{N}\C\right) \otimes S(\e_{\1}^{*})^{N}.
\end{align*}
This proves the desired result.
\end{proof}

Finally consider the case of the group $N'=\operatorname{Norm}_{G_{\0}}(\c_{\1})$ acting on the Cartan subspace $\c_{\1}.$  Note that since $x_{0} \in \c_{\1} \subseteq \e_{\1}$ one has that $\operatorname{Stab}_{G_{\0}}(\c_{\1})=H.$

\begin{theorem}\label{T:KRforc} Let $N'=\operatorname{Norm}_{G_{\0}}(\c_{\1}).$  Then there is a graded $N'$-module $\mathcal{H}\subseteq  S(\c_{\1}^{*})$ such that the multiplication map 
\[
\mathcal{H} \otimes S(\c_{\1}^{*})^{N'} \to S(\c_{\1}^{*})
\] gives an isomorphism of graded right $S(\c_{\1}^{*})^{N'}$-modules.  Furthermore, $\mathcal{H} \cong \Ind_{H}^{N'} \C$ as $N'$-modules.
\end{theorem}

\begin{proof}  Since $W$ is a pseudoreflection group one has that the multiplication map,
\[
\mathcal{H} \otimes S(\c_{\1}^{*})^{W} \to S(\c_{\1}^{*}),
\] gives an isomorphism as graded right $ S(\c_{\1}^{*})^{W}$-modules and  where $\mathcal{H}=\Ind_{1}^{W} \C$ \cite[Theorem 12.4.7]{goodman}.   However, viewing $\mathcal{H}$ as an $N'$-module by inflation through the canonical quotient map $N' \to N'/H=W,$ we have $\mathcal{H}=\Ind_{H}^{N'} \C$ as $N'$-modules.  
\end{proof}

\subsection{} \Label{SS:panyushev}\setcounter{Df}{0} Taken together Theorems~\THKRforg, 
\THKRfore, and ~\THKRforc\ imply that when $\g$ is both stable and polar the injective graded maps given in 
Theorems~\THlunarichardson(c)  and \THdadokkac(c)  are of free $R$-modules of the
same rank (namely, rank $\dim M^{H}$), where 
\[
R:=S(\g_{\1}^{*})^{G_{\0}}\cong S(\e_{\1}^{*})^{N}\cong S(\c_{\1}^{*})^{W}.
\]
In general such maps are proper injections but become isomorphisms over the field of fractions of $R.$  
However, one can apply a result of Panyushev \cite[Theorem 1]{panyushev} on covariants to show that in fact one need
only invert a single element of $R$. Arguing just as in the proof of \cite[Proposition 4]{panyushev}, one has the following result.

\begin{theorem} \Label{T:invertingD} Let $\g$ be stable and polar and let $M$ be a $G_{\0}$-module.  Let $D' \in R$ be such that for $y \in \g_{\1}$ $D'(y) \neq 0$ implies $y$ is generic.  Then the graded $R$-module maps induced by restriction 
\[
\rho: \Hom_{G_{\0}}(S(\g_{\1}),M) \to \Hom_{N}(S(\e_{\1}),M) \to  \Hom_{N'}(S(\c_{\1}),M)
\] are isomorphisms after one localizes at $D'.$  
\end{theorem}

Recall from Section~\ref{SS:dadokkac} that, for the classical Lie superalgebras which are stable and polar considered in Sections~\ref{SS:typeA}--\ref{SS:typeF4}, the set of generic points of $\g_{\1}$ is precisely the principal open set defined by the discriminant of $W.$  Therefore, in these cases, the above theorem applies when one localizes at the discriminant.

\section{Construction of Detecting Subalgebras}\Label{S:morecohom}

\subsection{} \Label{SS:detecting} \setcounter{Df}{0}\emph{For the remainder of the paper we assume $\g$ is stable and polar.}  Some results require this assumption while others still hold under weaker hypotheses.  To avoid technicalities and for clarity of exposition we assume both.   

We are now prepared to apply the invariant theory 
results of the previous section to our 
study of cohomology. In order to accomplish this goal we first construct two subsuperalgebras 
of $\g$ which will detect the relative cohomology ring of $\g$.  This is done as follows.
%Recall that if $M_{1}$ and $M_{2}$ are $K$-modules, then
%\begin{equation}\Label{E:homiso}
%\Hom_{\operatorname{Lie}(K)}(M_{1},M_{2})=\Hom_{K^{0}}(M_{1},M_{2}).
%\end{equation} 
  Since $\g$ is stable one can fix a generic point 
$x_{0} \in\g_{\1},$ and we let $H=G_{\0,x_{0}}$ and $N=\operatorname{Norm}_{G_{\0}}(H).$   

Let $\e_{\1}=\g_{\1}^{H}$ as in the previous section.  Let
\[
\e_{\0} = \operatorname{Lie}(N)\subseteq \g_{\0}
\]     and set 
\begin{equation}\Label{E:defe}
\e=\e_{\0} \oplus \e_{\1} \subseteq \g.
\end{equation}

Since $\g$ is polar, one can let $\c_{\1} \subseteq \g_{\1}$ be the Cartan subspace chosen 
in the previous section so that it satisfies \eqref{E:containments}.  Let
\[
\c_{\0} = \operatorname{Lie}(H) \subseteq \g_{\0 }.
\] Set 
\begin{equation}\Label{E:defc}
\c=\c_{\0} \oplus \c_{\1} \subseteq \g.
\end{equation}  Recall (cf.\ \cite[24.3.3-6]{tauvel}) that
\[
\operatorname{Lie}(H)=\{y \in \g_{\0} \:\vert\:  [y,x_{0}]=0\}.
\]  

\begin{theorem}\Label{T:isomofcohoms} Let ${\mathfrak g}$ be a classical Lie superalgebra which is stable and polar. Then there exists classical Lie subsuperalgebras $\e$ and $\c$ as constructed 
above such that the restriction maps induces isomorphisms of graded superalgebras
\[
\operatorname{H}^{\bullet}(\g,\g_{\0}; \C) \cong \operatorname{H}^{\bullet}(\e,\e_{\0}; \C)^{N/N^{0}} \cong 
\operatorname{H}^{\bullet}(\c,\c_{\0}; \C)^{W},
\] where $N/N^{0}$ and $W$ are finite pseudoreflection groups.  
In particular, $\operatorname{H}^{\bullet}(\g,\g_{\0}; \C)$ is 
isomorphic to a polynomial ring in $r:=\operatorname{dim}(\c_{\1})$ variables.  Additionally, each cohomology ring is integral over the image(s) of $res$ contained within it and so all three rings have the same Krull dimension.  
\end{theorem}

\begin{proof} First we prove that if $y,z \in \e_{\1},$ then $[y,z] \in \operatorname{Lie}(H)$. 
Note that by definition we have $x_{0} \in \e_{\1}$ and $[x_{0},x_{0}] \in \operatorname{Lie}(H).$  For $y \in \e_{\1}$ the super version of the Jacobi identity implies
\[
0=[y,[x_{0},x_{0}]]=[[y,x_{0}],x_{0}]-[x_{0},[y,x_{0}]]=2[[y,x_{0}],x_{0}].
\]  Thus, $[y,x_{0}] \in \operatorname{Lie}(H).$  For the general case let $y,z \in \e_{\1},$ then we have 
\[
[[y,z],x_{0}]=[y,[z,x_{0}]]-[[y,x_{0}],z]=0,
\] hence $[y,z] \in \operatorname{Lie}(H).$

From the statement above it can be verified that $\e$ and $\c$ are Lie subsuperalgebras of $\g.$
Since $N$ and $H$ are reductive subgroups of $G_{\0}$ it follows that $\e$ and $\c$ are classical Lie superalgebras.  
Consequently we can apply Theorem~\THcohomring\ to the pairs $(\e,\e_{\0})$ and $(\c,\c_{\0})$ and reinterpret 
Theorems~\THlunarichardson\ and \THdadokkac\  to obtain the stated results.
\end{proof} 

\subsection{Injectivity of cochains} \setcounter{Df}{0} We can also consider the relationship 
between relative cohomology for $\g,$ $\e,$ and $\c$ with coefficients in a finite dimensional supermodule other than the trivial supermodule.  

\begin{theorem} \Label{T:injectivityoncohains} Let $M$ be a finite dimensional $\g$-supermodule.  Then restriction induces injective graded maps
\begin{equation}\Label{E:injectivityoncochains}
C^{d}(\g,\g_{\0};M) \hookrightarrow C^{d}(\e,\e_{\0};M)^{N/N^{0}} \hookrightarrow C^{d}(\c,\c_{\0};M)^{W}.
\end{equation}
Furthermore, all three cochain complexes 
are free $\operatorname{H}^{\bullet}(\g,\g_{0};\C)$ modules of rank $\dim M^{H}.$
\end{theorem}

\begin{proof} First, recall that if $K$ is an algebraic group and $M_{1}$ and $M_{2}$ are $K$-modules, then $\Hom_{\operatorname{Lie}(K)}(M_{1},M_{2})=\Hom_{K^{0}}(M_{1},M_{2}).$  Using this to reinterpret \THlunarichardson(c) and \THdadokkac(c)\,
one has the first result.  
The second statement follows from Theorems~\THKRforg, \THKRfore, and \THKRforc.  
\end{proof}

\subsection{} \setcounter{Df}{0}In general the maps in Theorem~\ref{T:injectivityoncohains} fail to be surjective.  
However, more can be said in certain cases. 

\begin{theorem}\Label{T:krulldimoneiso}  Let $\g$ be a classical Lie superalgebra which is stable and polar.  Let $M$ be a 
finite dimensional $\g$-supermodule.  Furthermore, assume $\operatorname{H}^{\bullet}(\g,\g_{\0};\C)$ is isomorphic to a polynomial ring 
in one variable.  
Then restriction induces an isomorphism 
\[
\operatorname{H}^{d}(\g,\g_{\0};M) \cong \operatorname{H}^{d}(\e,\e_{\0};M)^{N/N^{0}}\cong \operatorname{H}^{d}(\c,\c_{\0};M)^{W}
\] for $d\gg 0.$

\end{theorem}

\begin{proof} First observe that if $R=\C[T]$ is a graded polynomial ring in one variable with 
$\operatorname{deg}(T)=t$, and 
$M \hookrightarrow N$ is an injective graded $R$-module homomorphism between 
two free $\Z$-graded $R$-modules of the same finite rank, then   
for $d\gg 0$ the map defines an isomorphism 
$M_{d} \cong N_{d}$  where $X_{d}$ denotes the $d$th graded piece of the module $X$. 
To see this one can compare dimensions.  Suppose that $U$ is a 
free graded $R$-module of rank $k:$
\[
U=Ru_{1} \oplus \dotsb \oplus Ru_{k}.
\]  Observe that for $d>\text{max}\{\operatorname{deg}(u_{1}), \dotsc, 
\operatorname{deg}(u_{k})\}$ one has 
\[
\operatorname{dim}\left(U_{d}\oplus \dotsb \oplus U_{d+t-1} \right)=k.
\]  
This along with the fact that $\operatorname{dim}(M_{d}) \leq \operatorname{dim}(N_{d})$ 
implies the claim.

Let us consider the first isomorphism in the statement of the theorem 
and leave the other nearly identical case to the reader.  Combining Theorem~\THinjectivityoncohains\ 
and the above observation one has
\[
C^{d}(\g,\g_{\0};M) \cong C^{d}(\e,\e_{\0};M)^{N/N^{0}}
\] for $d\gg 0.$  The result then follows from the fact that the differentials in \eqref{E:differential} are $N/N^{0}$ homomorphisms and, since $N/N^{0}$ is a finite group, taking $N/N^{0}$ invariants 
is an exact functor. 
\end{proof}

As seen from Table~1, the above stability result holds for $\mathfrak{gl}(1|n),$ $\mathfrak{q}(1),$ $D(2,1;\alpha),$ $G(3),$ and $F(4),$ among others.

\subsection{Localization of Cohomology}\Label{SS:isosafterlocalization}\setcounter{Df}{0} 
Let $D'$ be the element in $\operatorname{H}^{\bullet}(\g,\g_{\0};\C)$ chosen in Theorem~\ref{T:invertingD}.  If $X$ is a  $\operatorname{H}^{\bullet}(\g,\g_{\0};\C)$-module, write $X_{D'}$ to denote the localization of $X$ at $D'$.  
The following theorem relates the various relative cohomology groups under localization.

\begin{theorem} \Label{T:localizingcohom} Let $M$ be a finite dimensional $\g$-supermodule.  
Then restriction induces isomorphisms 
\[
\operatorname{H}^{\bullet}(\g,\g_{\0}; M)_{D'} \cong \left(\operatorname{H}^{\bullet}
(\e,\e_{\0}; M)^{N/N^{0}}\right)_{D'} \cong 
\left(\operatorname{H}^{\bullet}(\c,\c_{\0}; M)^{W}\right)_{D'}.
\]
\end{theorem}

\begin{proof}  Applying Theorem~\THinvertingD\ to \eqref{E:injectivityoncochains}, we see that 
\[
C^{d}(\g,\g_{\0};M)_{D'} \cong \left(C^{d}(\e,\e_{\0};M)^{N/N^{0}}\right)_{D'} \cong  \left(C^{d}(\c,\c_{\0};M)^{W}\right)_{D'}.
\]  The differentials are graded $\operatorname{H}^{\bullet}(\g,\g_{\0};\C)$-module homomorphisms 
and localization is an exact functor, thus one obtains the stated result.
\end{proof}

\section{Representation Theory for Detecting Subalgebras} 
\vskip .25cm 
In the previous section we constructed a Lie subsuperalgebra ${\mathfrak e}$ of ${\mathfrak g}$ 
which detects cohomology. This subalgebra will play an important role in defining rank varieties 
for finite dimensional supermodules. The definition of the rank varieties 
depends on understanding the representation theory of $\c.$  We develop the theory at this point. 

\subsection{The category $\mathcal{F}_{(\c,\c_{\0})}$}\Label{SS:categoryFL}\setcounter{Df}{0} In 
this subsection we investigate the general representation theory of the Lie superalgebra $\c$ defined in Section~\ref{SS:detecting}.   
The key fact which will be used below without comment is that 
\[
[\c_{\0},\c_{\1}]=0;
\] that is, elements of $\c_{\0}$ and $\c_{\1}$ commute.

\begin{lemma} \Label{L:dualizingprojs}  Let $M$ be a finite dimensional $\c_{\0}$-supermodule.  Then 
\begin{equation}\Label{E:dualiso}
\left( U(\c) \otimes_{U(\c_{\0})} M\right)^{*} \cong U(\c) \otimes_{U(\c_{\0})} M^{*}
\end{equation}
as $\c$-supermodules.  In particular, $\c$-supermodules are projective if and only if they are injective.
\end{lemma}

\begin{proof}  To prove \eqref{E:dualiso} we first consider the special case when $M=\C,$ the trivial supermodule.  
Fix a basis $x_{1},\dotsc, x_{r}$ for $\c_{\1}.$
By the PBW theorem for Lie superalgebras, if $N:=U(\c) \otimes_{U(\c_{\0})} \C,$ then the elements 
\[
x_{1}^{d_{1}}\dotsb x_{r}^{d_{r}} \otimes 1,
\] where $d_{1},\dotsc ,d_{r} \in \{0,1 \},$ form a basis for $N.$  Define a dual basis for $N^{*}$ by the formula
\[
\left[ x_{1}^{d_{1}}\dotsb x_{r}^{d_{r}} \otimes 1\right]^{*}\left( x_{1}^{e_{1}}\dotsb x_{r}^{e_{r}} \otimes 
1\right)=\delta_{d_{1},e_{1}}\dotsb \delta_{d_{r},e_{r}}.
\]  

Define the $\c$-supermodule homomorphism 
\[
\alpha: U(\c) \otimes_{U(\c_{\0})} \C \to \left( U(\c) \otimes_{U(\c_{\0})} \C\right)^{*},
\] to be the map induced by Frobenius reciprocity from the $\c_{\0}$-supermodule homomorphism $\C \to \left( U(\c) 
\otimes_{U(\c_{\0})} \C\right)^{*}$ given by $1 \mapsto \left[ x_{1}\dotsb x_{r} \otimes 1\right]^{*}.$  

We will prove $\alpha$ is an isomorphism.  Before doing so, let us set some notation.  Given a tuple 
$d=(d_{1}, \dotsc ,d_{r}) \in \{0,1 \}^{r},$ let $|d|=d_{1}+\dotsb +d_{r}$ and let $\tilde{d}_{i}=1-d_{i}.$  We claim 
\begin{equation}\Label{E:alphamap}
\alpha\left( x_{1}^{d_{1}}\dotsb x_{r}^{d_{r}} \otimes 1\right)=\pm \left[ x_{1}^{\tilde{d}_{1}}\dotsb 
x_{r}^{\tilde{d}_{r}} \otimes 1\right]^{*}.
\end{equation}  Note that if $\alpha$ satisfies \eqref{E:alphamap}, then it is an isomorphism. 

To prove \eqref{E:alphamap}, one inducts on $|d|,$ with the case $|d|=0$ being clear from the definition of $\alpha.$  
Now consider the case $|d|> 0.$  Fix $k$ so that $d_{k}=1$ and $d_{i}=0$ for $i=1,\dotsc ,k-1.$  Then,
\begin{align*}
\alpha\left( x_{1}^{d_{1}}\dotsb x_{r}^{d_{r}} \otimes 1\right)&=(-1)^{\p{\alpha}}x_{k}\alpha\left( x_{1}^{d_{1}}\dotsb x_{k-1}^{d_{k-1}}x_{k}^{\tilde{d}_{k}}x_{k+1}^{d_{k+1}}\dotsb x_{r}^{d_{r}} \otimes 1 \right)\\
     &=\pm x_{k}\left[ x_{1}^{\tilde{d}_{1}}\dotsb x_{k-1}^{\tilde{d}_{k-1}}x_{k}^{d_{k}}x_{k+1}^{\tilde{d}_{k+1}}\dotsb x_{r}^{\tilde{d}_{r}} \otimes 1\right]^{*},
\end{align*} where the second equality is by the inductive assumption.  Now on the one hand we have
\begin{align*}
\left( x_{k}\left[ x_{1}^{\tilde{d}_{1}}\dotsb x_{k}^{d_{k}}\dotsb  x_{r}^{\tilde{d}_{r}} \otimes 1\right]^{*}\right)&\left(  x_{1}^{e_{1}}\dotsb x_{r}^{e_{r}} \otimes 1\right) \\
&= \pm \left[ x_{1}^{\tilde{d}_{1}}\dotsb x_{k}^{d_{k}} \dotsb  x_{r}^{\tilde{d}_{r}} \otimes 1\right]^{*}\left(x_{k} x_{1}^{e_{1}}\dotsb x_{r}^{e_{r}} \otimes 1 \right)\\
 &= \pm \left[ x_{1}^{\tilde{d}_{1}}\dotsb  x_{k}^{d_{k}} \dotsb x_{r}^{\tilde{d}_{r}} \otimes 1\right]^{*}\left( x_{1}^{e_{1}}\dotsb x_{k-1}^{e_{k-1}}x_{k}^{e_{k}+1}x_{k+1}^{e_{k+1}}\dotsb  x_{r}^{e_{r}} \otimes 1 \right) \\
&=\begin{cases} \pm \;\delta_{\tilde{d}_{1},e_{1}}\dotsb \delta_{d_{k},e_{k}+1} \dotsb \delta_{\tilde{d}_{r},e_{r}}, & \text{if } e_{k}=0;\\
0, & \text{if } e_{k}=1.
\end{cases}
\end{align*}  Note that the second case follows from the observation that
\[
x_{k}^{2}=\frac{1}{2}[x_{k},x_{k}] \in \c_{\0}.
\]
On the other hand, 
\begin{equation*}
\left[ x_{1}^{\tilde{d}_{1}}\dotsb x_{r}^{\tilde{d}_{r}} \otimes 1\right]^{*}\left(  x_{1}^{e_{1}}\dotsb x_{r}^{e_{r}} \otimes 1\right) =
\delta_{\tilde{d}_{1},e_{1}}\dotsb \delta_{\tilde{d}_{k},e_{k}} \dotsb \delta_{\tilde{d}_{r},e_{r}}. 
\end{equation*}  Using the assumption that $d_{k}=1,$ one can verify that these two expressions are equal up to sign.  This proves that $\alpha$ defines an isomorphism.  The general case then follows by an argument using the tensor identity.

Now consider the second statement.  Let $P$ be a finite dimensional projective $\c$-supermodule.  Then by Frobenius reciprocity we have an even surjective map
\[
N:=U(\c)\otimes_{U(\c_{\0})} P \to P.
\]  Since $P$ is projective, we have an even isomorphism $N\cong P \oplus U.$  Dualizing and applying \eqref{E:dualiso} we see that $P^{*}$ is a direct summand of the supermodule $U(\c)\otimes_{U(\c_{\0})} P^{*}.$  However $U(\c)\otimes_{U(\c_{\0})} P^{*}$ is projective since $P^{*}$ is an object of $\mathcal{F}_{(\c,\c_{\0})}$, hence a projective $\c_{\0}$-supermodule, and induction is exact.  Therefore $P^{*}$ is both projective and injective, hence so is $P.$ 
\end{proof}

\begin{lemma}\Label{L:blocks}  Let $M$ be a finite dimensional $\c$-supermodule, and let 
\begin{equation}\Label{E:isotypy}
M=M_{1} \oplus \dotsb \oplus M_{t}
\end{equation}
be a decomposition into $\c_{\0}$-isotypic components.  
Then \eqref{E:isotypy} is a decomposition into $\c$-supermodules.
\end{lemma}

\begin{proof} Let $L$ be an simple $\c_{\0}$-supermodule which appears as a direct summand of the isotypic component $M_{k}$ of $M$.  Let $y \in \c_{\1}.$   Then there is a $\c_{\0}$-supermodule homomorphism $L \to M$ given by left multiplication by $y.$  By the irreducibility of $L$ the image of $L$ under this map is either zero or isomorphic to $L.$  In either case, $yL \subseteq M_{k}.$  Therefore, for any $y \in \c_{\1}$ we have $yM_{k} \subseteq M_{k}.$   This, along with the observation that the decomposition \eqref{E:isotypy} respects the $\Z_2$-grading of $M$, implies the decomposition \eqref{E:isotypy} is as $\c$-supermodules. 
\end{proof}

Suppose an $\c$-supermodule $M$ contains one and only one simple $\c_{\0}$-supermodule as a composition factor, say $L.$  In this case we say \emph{$M$ is of type L.}  Let us write $\mathcal{F}_{L}$ for the full subcategory of $\mathcal{F}_{(\c,\c_{\0})}$ consisting of all $\c$-supermodules of type $L.$   Applying Lemma~\ref{L:blocks} we have the following direct sum decomposition of the category $\mathcal{F}_{(\c,\c_{\0})},$
\begin{equation}\Label{E:blockdecomp}
\mathcal{F}_{(\c,\c_{\0})}=\bigoplus_{L} \mathcal{F}_{L},
\end{equation}
where the sum is over a complete, irredundant set of simple $\c_{\0}$-supermodules.  We refer to the subcategory $\mathcal{F}_{L}$ as the \emph{block of $\mathcal{F}_{(\c,\c_{\0})}$ of type $L.$}

\begin{lemma}\Label{L:irrepsinFL} Let $L$ be a finite dimensional simple $\c_{\0}$-supermodule.  The finite dimensional simple $\c$-supermodules which lie in $\mathcal{F}_{L}$ are precisely the composition factors of 
\[
U(\c) \otimes_{U(\c_{\0})} L.
\]  In particular, the block $\mathcal{F}_{L}$ contains only finitely many simples.  
\end{lemma}

\begin{proof}  Let $T$ be an simple $\c$-supermodule lying in $\mathcal{F}_{L}.$  Then by Frobenius reciprocity we have a surjective map
\[
\bigoplus U(\c) \otimes_{U(\c_{\0})} L \to T,
\] where the sum runs over some finite index set.  Thus $T$ is a composition factor of the supermodule $U(\c) \otimes_{U(\c_{\0})} L.$
\end{proof}

\subsection{Representations for rank one Lie superalgebras} \Label{SS:rankone}\setcounter{Df}{0} 
The definition of rank varieties in Section~\ref{SS:rankvariety} depends on restricting to 
rank one Lie superalgebras. Here we compile some basic results on finite dimensional representations of these superalgebras.  
Let $\x$ be a Lie superalgebra generated by a single odd vector, $x.$  There are the 
following two possibilities for $\x.$

\begin{enumerate}
\item [I.]  We have $[x,x]=0.$  Then $\dim \x_{\0}=0,$ $\dim \x_{\1}=1,$ and $\x$ is a one dimensional abelian Lie superalgebra concentrated in degree $\1.$
\item [II.]  We have $h:=[x,x] \neq 0.$  Then $\dim \x_{\0}=1,$ $\dim \x_{\1} =1,$ and $\x$ is isomorphic to the Lie superalgebra $\mathfrak{q}(1).$
\end{enumerate}

Let us first consider Case~I.  The universal enveloping superalgebra $U(\x)$ is spanned by $1$ and $x$ and subject to the relation $x^{2}=0.$  It 
is easy to see that the trivial module is the only finite dimensional simple $U(\x)$-supermodule, and $U(\x)$ is the projective cover of $\C.$  
Furthermore, a direct calculation verifies that $U(\x)^{*} \cong U(\x)$ as $U(\x)$-supermodules, so $U(\x)$ is also injective.  From this we see that $U(\x)$ is the unique self extension of $\C$ and that $\left\{U(\x),\C \right\}$ is a complete, irredundant list of 
indecomposable $U(\x)$-supermodules.

Now consider Case~II.  By the PBW theorem for Lie superalgebras 
$U(\x)$ is spanned by monomials of the form 
\[
x^{r}h^{s}
\] where $r \in \{0,1 \}$ and $s \in \mathbb{Z}_{\geq 0}.$  In particular, one has
an induction functor from the category of finite dimensional $\x_{\0}$-supermodules to the category of 
finite dimensional $\x$-supermodules given by 
\[
\Ind_{\x_{\0}}^{\x} M = U(\x) \otimes_{U(\x_{\0})} M.
\]  Furthermore, this functor is both left adjoint to restriction and exact because 
$U(\x)$ is a free right $U(\x_{\0})$-supermodule.

 From the PBW theorem one sees that $U(\x_{\0})$ is isomorphic to the polynomial ring in one variable.  Consequently, the category of finite dimensional $U(\x_{\0})$-supermodules is semisimple and we obtain the following classification of the simple $\x_{\0}$-supermodules:  for each $\lambda \in \C$ let $\C_{\lambda}$ denote $\C$ viewed as an $\x_{\0}$-supermodule concentrated in degree $\0$ where $h$ acts by the scalar $\lambda,$ then the set $\{ \C_{\lambda}\:\vert\:  \lambda \in \C \}$ is a complete irredundant set of simple $\x_{\0}$-supermodules.  Furthermore, since $U(\x_{\0})$ is semisimple, these are projective as well.  Since induction 
is exact the $\x$-supermodule
\begin{equation}\Label{E:Plambda}
P(\lambda):=\Ind_{\x_{0}}^{\x} \C_{\lambda} = U(\x) \otimes_{U(\x_{\0})} \C_{\lambda}
\end{equation}
is projective for each $\lambda \in\C.$

One the one hand, by Frobenius reciprocity every simple $U(\x)$-supermodule is the quotient of 
$P(\lambda)$ for some $\lambda$.  
On the other hand, a direct calculation verifies that $P(\lambda)$ is a simple 
supermodule if and only if $\lambda \neq 0.$  To warn the reader that there are subtleties, note that $P(\lambda)$ is the direct 
sum of two $1$-dimensional $U(\x)$-modules when $\lambda \neq 0$.   In the case when $\lambda=0$ we have that $P(\lambda)$ has simple head and socle, each isomorphic to the trivial supermodule.

Now consider indecomposables.  By the above discussion, all the simple supermodules 
lie in separate blocks.  Furthermore, since all simple supermodules except the trivial supermodule are projective, the only indecomposable in the non-principal blocks is the simple supermodule itself.

 Now consider the principal block.  Applying Theorem~\THcohomring\ we see that 
\[
\dim\left(\Ext_{\mathcal{F}_{(\x,\x_{\0})}}^{1}(\mathbb{C}, \mathbb{C}) \right)=1.
\]  That is, $P(0)$ is the unique (up to isomorphism) self extension of the trivial supermodule.  Consequently, $P(0)$ is self dual, hence injective.  
From this observation it follows that $P(0)$ and the trivial supermodule are the only indecomposable supermodules in the principal block.

We summarize the above analysis in the following proposition.  

\begin{prop} \Label{L:irrepsofQ(1)}  If $\x$ is as in Case~I above, then the following statements about the category of finite dimensional 
$\x$-supermodules hold.

\begin{itemize}
\item[(a)] The trivial supermodule $L(0):=\C$ is the only simple supermodule.
\item[(b)] The projective cover of $L(0)$ is $P(0):=U(\x).$
\item[(c)] The supermodule $P(0)$ is self dual, hence injective.
\item[(d)] The set $\left\{P(0),L(0) \right\}$ is a complete set of indecomposable supermodules.
\end{itemize}

If $\x$ is as in Case~II above, then the following statements about the category of finite dimensional $\x$-supermodules hold.  
Given $\lambda \in \C,$ let $P(\lambda)$ be as in \eqref{E:Plambda}.  Let $L(\lambda)$ denote the head of $P(\lambda).$
\begin{itemize}
\item[(e)] The set $\{L(\lambda)\:\vert\:  \lambda \in \C \}$ is a complete set of simple $\x$-supermodules.
\item[(f)] For all $\lambda \in \C,$ $P(\lambda)$ is the projective cover of $L(\lambda).$
\item[(g)] If $\lambda \neq  0,$ then $L(\lambda)=P(\lambda)$.
\item[(h)] For all $\lambda \in \mathbb{C},$ $P(\lambda)$ is dual to $P(-\lambda)$, hence injective.
\item[(i)] The set $\{L(\lambda)\:\vert\: \lambda \in \mathbb{C} \} \cup \{P(0) \}$ is a complete set of indecomposable supermodules.  
\end{itemize}

We remark that in both cases the supermodules $P(\lambda)$ always satisfy 
\[
\dim (P(\lambda)_{\0})=\dim (P(\lambda)_{\1})=1.
\]
In particular, $P(\lambda)$ is always two dimensional.
\end{prop}

\section{Support Varieties} \Label{S:supportvarieties}

\subsection{} \Label{SS:def}\setcounter{Df}{0} In this section we define the notion of the support variety  
of a finite dimensional $\g$-supermodule $M$ and study the properties of these varieties.
Let $\mathfrak{a}$ be a classical Lie superalgebra (e.g. $\mathfrak{a}$ could be 
one of $\g,$ $\e,$ or $\c$ from the previous sections).  Let $M$ and $N$ be $\mathfrak{a}$-supermodules 
in the category $\mathcal{F}_{(\mathfrak{a},\mathfrak{a}_{\0})}$.  Recall that by 
Theorem~\THfinitegen2\ one knows that

\[
\Ext_{\mathcal{F}_{(\mathfrak{a},\mathfrak{a}_{\0})}}^{\bullet}(M,N) \cong \operatorname{H}^{\bullet}(\mathfrak{a}, \mathfrak{a}_{\0}; M^{*} \otimes N)
\] 
 is a finitely generated $\operatorname{H}^{\bullet}(\mathfrak{a}, \mathfrak{a}_{\0};\C)$-module.  Let 
\[
I_{(\mathfrak{a},\mathfrak{a}_{\0})}(M,N)=\operatorname{Ann}_{\operatorname{H}^{\bullet}(\mathfrak{a},
\mathfrak{a}_{\0};\C)}(\operatorname{H}^{\bullet}(\mathfrak{a}, \mathfrak{a}_{\0}; M^{*} \otimes N)),
\] the annihilator ideal of this module.  We define the \emph{relative support variety of the pair $(M,N)$} to be 
\begin{equation}\Label{E:suppvardef} \mathcal{V}_{(\mathfrak{a},\mathfrak{a}_{\0})}(M,N)=
\operatorname{MaxSpec}(\operatorname{H}^{\bullet}(\mathfrak{a}, \mathfrak{a}_{\0};\C)/I_{(\mathfrak{a},\mathfrak{a}_{\0})}(M,N)),
\end{equation} the maximal ideal spectrum of the quotient of $\operatorname{H}^{\bullet}(\mathfrak{a}, \mathfrak{a}_{\0};\C)$ by $I_{(\mathfrak{a},\mathfrak{a}_{\0})}(M,N).$
 For short when $M=N$, write 
\begin{align*}
I_{(\mathfrak{a},\mathfrak{a}_{\0})}(M)&=I_{(\mathfrak{a},\mathfrak{a}_{\0})}(M,M),\\
\mathcal{V}_{(\mathfrak{a},\mathfrak{a}_{\0})}(M)&=\mathcal{V}_{(\mathfrak{a},\mathfrak{a}_{\0})}(M,M).
\end{align*}  We call $\mathcal{V}_{(\mathfrak{a},\mathfrak{a}_{\0})}(M)$ the \emph{support variety} of $M.$  
Let us remark that, just as for finite groups, $I_{(\mathfrak{a},\mathfrak{a}_{\0})}(M)$ is precisely 
the annihilator ideal of the identity element of $\operatorname{H}^{\bullet}(\mathfrak{a},\mathfrak{a}_{\0};M^{*}\otimes M)$ 
viewed as a ring under the Yoneda product.  

 Set $r:=\operatorname{dim}(\c_{\1}).$ If $\mathfrak{a}$ equals $\g,$ $\e,$ or $\c,$ then $\mathcal{V}_{(\mathfrak{a},\mathfrak{a}_{\0})}(M) \subseteq \mathcal{V}_{(\mathfrak{a},\mathfrak{a}_{\0})}(\C) \cong \mathbb{A}^{r}.$   That is, $\mathcal{V}_{(\mathfrak{a},\mathfrak{a}_{\0})}(M)$ can naturally 
be viewed as a conical (since the defining ideal is graded) affine subvariety of $\mathbb{A}^{r}.$ 

Recall that when $\g$ is stable and polar one has canonical restriction maps 
\[
\operatorname{H}^{\bullet}(\g,\g_{\0}; \C) \to \operatorname{H}^{\bullet}(\e,\e_{\0}; \C) \to \operatorname{H}^{\bullet}(\c,\c_{\0}; \C)
\] which induce maps, which we call $\res^{*},$
\[
\mathcal{V}_{(\c,\c_{\0})}(\C) \to \mathcal{V}_{(\e,\e_{\0})}(\C) \to \mathcal{V}_{(\g,\g_{\0})}(\C).
\] 
By Theorem~\ref{T:isomofcohoms} one then has 
\begin{equation}\Label{E:resiso}
\mathbb{A}^{r}\cong \mathcal{V}_{(\c,\c_{\0})}(\C)/W \cong \mathcal{V}_{(\e,\e_{\0})}(\C)/(N/N^{0}) \cong \mathcal{V}_{(\g,\g_{\0})}(\C).
\end{equation}

Furthermore, if $M$ is a finite dimensional $\g$-supermodule, then $\res^{*}$ restricts to give maps,
\[
\mathcal{V}_{(\c,\c_{\0})}(M) \to \mathcal{V}_{(\e,\e_{\0})}(M) \to \mathcal{V}_{(\g,\g_{\0})}(M).
\] Since the first two varieties are stable under the action of $W$ and $N/N^{0},$ respectively, 
applying \eqref{E:resiso} we have the following embeddings of varieties induced by $\res^{*}$,
\begin{align}\Label{E:resmapsbetweenvarieties}
\mathcal{V}_{(\c,\c_{\0})}(M)/W &\cong \res^{*}\left(\mathcal{V}_{(\c,\c_{\0})}(M)\right)\subseteq  \mathcal{V}_{(\g,\g_{\0})}(M), \\
\mathcal{V}_{(\e,\e_{\0})}(M)/(N/N^{0}) &\cong \res^{*}\left(\mathcal{V}_{(\e,\e_{\0})}(M)\right)\subseteq \mathcal{V}_{(\g,\g_{\0})}(M). \notag
\end{align}

\subsection{Relating support varieties}\Label{SS:relatingvarieties}\setcounter{Df}{0} Naturally one would like to better understand 
the relationship between the varieties $ \mathcal{V}_{(\g,\g_{\0})}(M)$, $ \mathcal{V}_{(\e,\e_{\0})}(M),$ and 
$ \mathcal{V}_{(\c,\c_{\0})}(M)$ for a finite dimensional $\g$-supermodule $M.$  In particular, one would like to understand 
the maps given in \eqref{E:resmapsbetweenvarieties}.  It may be that these maps are surjective, at least for suitably nice supermodules.  In any case, one can use results from Section~\Smorecohom\ to gain some insight. 

\begin{theorem} \label{T:rankonecase} Assume $\g$ is a classical Lie superalgebra which is stable and polar.  Assume further that $\operatorname{H}^{\bullet}(\g,\g_{\0};\C)$ is isomorphic to a polynomial ring in one variable.  
Let $M$ be a finite dimensional $\g$-supermodule.  Then $\res^{*}$ defines the following isomorphisms
\[
\res^{*}\left( \mathcal{V}_{(\c,\c_{\0})}(M)\right) = \res^{*}\left(\mathcal{V}_{(\e,\e_{\0})}(M)\right) = \mathcal{V}_{(\g,\g_{\0})}(M).
\]
\end{theorem}

\begin{proof}  We consider the map $\res^{*}:\mathcal{V}_{(\c,\c_{\0})}(M) \to \mathcal{V}_{(\g,\g_{\0})}(M).$  
The other case follows similarily.  

In light of \eqref{E:resmapsbetweenvarieties}, it suffices to prove $\res^{*}$ is surjective.  
Also note that by the ``lying over'' theorem for maximal ideals \cite[Theorem 5.10, Corollary 5.8]{atiyah}, 
$\res^{*}(\mathcal{V}_{(\c,\c_{\0})}(M))$ is precisely the subvariety of $\mathcal{V}_{(\g,\g_{\0})}(\C)$ defined by the ideal $\res^{-1}(I_{(\c,\c_{\0})}(M)).$  
We claim that 
\begin{equation}\Label{E:idealinclusions}
I_{(\g,\g_{\0})}(M) \subseteq \res^{-1}(I_{(\c,\c_{\0})}(M)) \subseteq \sqrt{I_{(\g,\g_{\0})}(M)},
\end{equation}
where the latter is the radical ideal of $I_{(\g,\g_{\0})}(M).$  This suffices to prove the desired 
assertion as the radical of an ideal defines the same variety as the ideal itself.

First, by Theorem~\THkrulldimoneiso, one can fix $N>0$ so that the restriction map 
\[
\operatorname{H}^{d}(\g,\g_{\0};M) \to  \operatorname{H}^{d}(\c,\c_{\0};M)^{W} \subseteq \operatorname{H}^{d}(\c,\c_{\0};M)
\] is an isomorphism for all $d\geq N.$  

The first inclusion of \eqref{E:idealinclusions} is clear from the description of $I_{(\g,\g_{\0})}(M)$ as the annihilator of the identity of $\operatorname{H}^{\bullet}(\g,\g_{\0};M^{*}\otimes M)$.  On the other hand, let $x \in \res^{-1}(I_{(\c,\c_{\0})}(M))$ and 
let $m \in \operatorname{H}^{p}(\g; \g_{\0}; M^{*} \otimes M)$ for some $p \geq 0.$  Then $\res(x)\res(m)=0$ so 
\[
\res(x^{N}m)=\res(x)^{N}\res(m)=0.
\]
However, by our choice of $N$ the map $\res$ is injective, so $x^{N}m=0.$  But since $m$ was arbitrary, 
we see that $x^{N} \in I_{(\g,\g_{\0})}(M);$ that is, $x \in \sqrt{I_{(\g,\g_{\0})}(M)}.$
\end{proof}

Before continuing, we recall from commutative algebra (cf.\  
\cite[Chapter 3, Exercise 1]{atiyah}) that if $R$ is a commutative ring, $S \subseteq R$ a 
multiplicatively closed set, and $U$ a finitely generated $R$-module, then $S^{-1}U=0$ if and only if $sU=0$ for 
some $s \in S$.   
In particular, this implies that if $\mathfrak{m} \in \mathcal{V}_{(\g,\g_{\0})}(\C),$ 
then  $\mathfrak{m} \in \mathcal{V}_{(\g,\g_{\0})}(M)$ if and only if 
$\operatorname{H}^{\bullet}(\g,\g_{\0};M^{*}\otimes M)_{\mathfrak{m}}\neq 0.$

In general, \eqref{E:resmapsbetweenvarieties} is ``generically'' an isomorphism in the 
following sense.  Let $D' \in \operatorname{H}^{\bullet}(\g,\g_{\0};\C)$ be as chosen in Theorem~\ref{T:invertingD}
and let 
\[
\mathcal{O}_{D'}=\setof{a \in \mathcal{V}_{(\g,\g_{\0})}(\C)}{D'(a) \neq 0}=\setof{\mathfrak{m}\in \operatorname{MaxSpec}(\operatorname{H}^{\bullet}(\g,\g_{\0};\C))}{D' \not\in \mathfrak{m}},
\] the principal dense open set defined by $D'.$  In particular, as remarked in Section~\ref{SS:dadokkac}, for the classical Lie superalgebras considered in Sections~\ref{SS:typeA}--\ref{SS:typeF4} which are stable and polar one can take $D'$ to be the discriminant of $W$ and then $\mathcal{O}_{D'}$ is precisely the image in $\mathcal{V}_{(\g,\g_{\0})}(\C)=\g_{\1}/G_{\0}$ of the set of generic points of $\g_{\1}$.

\begin{theorem}\Label{T:genericequality} Let $M$ be a finite dimensional $\g$-supermodule.  Then 
\[
\res^{*}\left( \mathcal{V}_{(\c,\c_{\0})}(M)\right) \cap \mathcal{O}_{D'} = 
\res^{*}\left(\mathcal{V}_{(\e,\e_{\0})}(M)\right) \cap \mathcal{O}_{D'} = \mathcal{V}_{(\g,\g_{\0})}(M)\cap \mathcal{O}_{D'}.
\]
\end{theorem}

\begin{proof}  We prove $\res^{*}\left( \mathcal{V}_{(\c,\c_{\0})}(M)\right) \cap \mathcal{O}_{D'} = 
\mathcal{V}_{(\g,\g_{\0})}(M)\cap \mathcal{O}_{D'}$ and leave the other case for the reader.

It is clear that $\res^{*}\left( \mathcal{V}_{(\c,\c_{\0})}(M)\right) 
\cap \mathcal{O}_{D'} \subseteq 
\mathcal{V}_{(\g,\g_{\0})}(M)\cap \mathcal{O}_{D'}.$  On the other hand, let $\mathfrak{m} \in 
\mathcal{V}_{(\g,\g_{\0})}(M)\cap \mathcal{O}_{D'}.$  That is, 
$\operatorname{H}^{\bullet}(\g,\g_{\0};M^{*}\otimes M)_{\mathfrak{m}}\neq 0,$ and 
$D' \not\in \mathfrak{m}.$  Theorem~\THlocalizingcohom\ implies that 
\[
0 \neq \operatorname{H}^{\bullet}(\g,\g_{\0};M^{*}\otimes M)_{\mathfrak{m}} \cong 
\left( \operatorname{H}^{\bullet}(\c,\c_{\0};M^{*}\otimes M)^{W}\right)_{\mathfrak{m}} \subseteq \operatorname{H}^{\bullet}(\c,\c_{\0};M^{*}\otimes M)_{\mathfrak{m}}.
\] This implies the reverse inclusion.
\end{proof}

\subsection{Rank Varieties} \Label{SS:rankvariety} \setcounter{Df}{0} Given Theorems~\ref{T:rankonecase} and~\ref{T:genericequality}, we are motivated to more closely study the support varieties associated to the detecting subsuperalgebra $\c.$   The 
goal of this subsection is to define rank varieties and prove an analogue the 
Avrunin-Scott Theorem (first conjectured for finite groups by Carlson).  Our approach is based on  
ideas used in the context of restricted Lie algebras by Friedlander and Parshall \cite[Theorem 2.7]{FP}.

Let $M$ be a finite dimensional $\c$-supermodule.  As a matter of notation, given a homogeneous 
element $x \in \c,$ let $\langle x \rangle$ denote the Lie subsuperalgebra generated by $x.$   
Define the \emph{rank variety} of $M$ to be 
\[
{\mathcal V}_{\c}^{\text{rank}}(M)=\left\{ x \in\c_{\1}\:\vert\:  M
\text{ is not projective as an $\langle x \rangle$-supermodule} \right\} \cup \{0 \}.
\]
Observe that 
\[
{\mathcal V}_{\c}^{\text{rank}}(\mathbb{C})=\c_{\1}.
\]

We record some basic properties of rank varietes in the following proposition.
\begin{prop}\label{L:tensorproduct} Let $M$ and $N$ be finite dimensional $\c$-supermodules.  Then,
\begin{itemize}
\item[(a)] ${\mathcal V}_{\c}^{\operatorname{rank}}(M \otimes N)=
{\mathcal V}_{\c}^{\operatorname{rank}}(M) \cap {\mathcal V}_{\c}^{\operatorname{rank}}(N).$ 
\item[(b)] ${\mathcal V}_{\c}^{\operatorname{rank}}(M^{*})={\mathcal V}_{\c}^{\operatorname{rank}}(M)$
\item[(c)] ${\mathcal V}_{\c}^{\operatorname{rank}}(M^{*} \otimes M)={\mathcal V}_{\c}^{\operatorname{rank}}(M)$.
\end{itemize}
\end{prop}
\begin{proof}  (a) From Proposition~\LirrepsofQ(1)\ one observes that $x \in {\mathcal V}^{\operatorname{rank}}_{\c}(M)$ if and only 
if $M$ contains a trivial direct summand as a $\langle x \rangle$-supermodule. A direct calculation verifies 
that if $P(\lambda)$ and $P(\mu)$ are two indecomposable projective $\langle x \rangle$-supermodules, then 
\[
P(\lambda) \otimes P(\mu) \cong P(\lambda+\mu) \oplus P(\lambda + \mu)
\] as $\langle x \rangle$-supermodules.  The result follows from these observations.

(b)--(c) From Proposition~\LirrepsofQ(1)\ we see that for 
any $x \in \c_{\1},$ $M^{*}$ is projective as an $\langle x \rangle$-supermodule if and only if $M$ is projective.  This implies the first result.  The second result then follows from (a). 
\end{proof}

\begin{theorem} \Label{T:avruninscott} Let $M$ be a finite dimensional $\c$-supermodule.  Then there is an isomorphism 
\[ {\mathcal V}_{\c}^{\operatorname{rank}}(\mathbb{C}) \to {\mathcal V}_{(\c, \c_{\0})}(\mathbb{C}) 
\] which restricts to give an isomorphism 
\[
{\mathcal V}_{\c}^{\operatorname{rank}}(M)\to {\mathcal V}_{(\c, \c_{\0})}(M)   
\] for any $\c$-supermodule $M$.
\end{theorem}

\begin{proof}  The isomorphism ${\mathcal V}_{(\c, \c_{\0})}(\mathbb{C}) \cong 
{\mathcal V}_{\c}^{\text{rank}}(\mathbb{C})$ is a direct consequence of  Hilbert's Nullstellensatz as
$\operatorname{H}^{\bullet}(\c,\c_{\0};\C)\cong S(\c_{\1}^{*}).$

Now let $M$ be a finite dimensional $\c$-supermodule.  To prove the second statement, we identify ${\mathcal V}_{(\c, \c_{\0})}(\mathbb{C}) $ with $\c_{\1}$ using the above isomorphism and show that $x \in {\mathcal V}_{(\c, \c_{\0})}(M)$ if and only if $x \in {\mathcal V}^{\operatorname{rank}}_{\c}(M).$ 

First, let $0 \neq x \in {\mathcal V}^{\operatorname{rank}}_{\c}(M).$
Using Proposition~\LirrepsofQ(1)\ one can verify directly that the isomorphism $ {\mathcal V}_{\langle x \rangle}^{\text{rank}}(\mathbb{C}) \to {\mathcal V}_{(\langle x \rangle, \langle x \rangle_{\0})}(\mathbb{C})$
restricts to give an isomorphism  
\begin{equation*} {\mathcal V}_{\langle x \rangle}^{\text{rank}}(M) \xrightarrow{\cong}
{\mathcal V}_{(\langle x \rangle, \langle x \rangle_{\0})}(M).\end{equation*}
By choosing coordinates one can verify 
that the restriction map 
\begin{equation*}
\res^{*}:{\mathcal V}_{(\x,\x_{\0})}(\mathbb{C})=\C x \to {\mathcal V}_{(\c,\c_{\0})}(\mathbb{C}) =\c_{\1}
\end{equation*}
is the injective map which has as its image the line which goes through $0$ and $x.$  By restricting this map one then has the injective map 
\begin{equation*}
\res^{*}:{\mathcal V}_{(\x, \x_{\0})}(M) \hookrightarrow {\mathcal V}_{(\c, \c_{\0})}(M)
\end{equation*}
with image equal to the line through $0$ and $x.$
Therefore, if $0 \neq x \in {\mathcal V}^{\operatorname{rank}}_{\c}(M),$ then  $x \in {\mathcal V}_{(\c, \c_{\0})}(M).$

Conversely, assume $x \not\in  {\mathcal V}_{\c}^{\text{rank}}(M).$  Then $M^{*}\otimes M$ is projective as 
an $\x$-supermodule and so $\operatorname{H}^{q}(\x,\x_{\0};M^{*}\otimes M)=0$ for $q>0$.  Let $I=\langle x \rangle + \c_{\0},$ a graded ideal of $\c.$  Viewing $I$ as a Lie subsuperalgebra of $\c$ and $M^{*} \otimes M$ as a $I$-supermodule by restriction, we wish to prove $\operatorname{H}^{q}(I,I_{\0};M^{*}\otimes M)=0$ for $q>0$.  

Observe that the cochains used to define $\operatorname{H}^{\bullet}(I,I_{\0};M^{*}\otimes M)$ in 
Section~\ref{SS:lierelcohom} can be obtained by taking the $\c_{\0}$ invariants of the cochains which define $\operatorname{H}^{\bullet}(\x,\x_{\0};M^{*}\otimes M)$.  Since the cochains which define $\operatorname{H}^{\bullet}(\x,\x_{\0};M^{*}\otimes M)$ are finitely semisimple as $\c_{\0}$-supermodules, taking $\c_{\0}$ invariants is an exact functor.  Furthermore, using that $[z,x]=0$ for any $z \in \c_{\0},$ one sees that the differentials are $\c_{\0}$-supermodule homomorphisms.  Taken together, this discussion implies $\operatorname{H}^{q}(I,I_{\0};M^{*}\otimes M) \cong \operatorname{H}^{q}(\x,\x_{\0};M^{*}\otimes M)^{\c_{\0}}=0$ for all $q > 0.$

There is a Hochschild-Serre spectral sequence with $E_{2}$ page
\[
E_{2}^{p,q}=\operatorname{H}^{p}(\c/I,0;\operatorname{H}^{q}(I,I_{\0};M^{*}\otimes M))\Rightarrow 
\operatorname{H}^{\bullet}(\c,\c_{\0};M^{*}\otimes M)
\]
(cf.\ \cite[Sections 1.5--1.6]{fuks}).  Since $\operatorname{H}^{q}(I,I_{\0};M^{*}\otimes M)=0$ for $q>0,$ one has 
\[
\operatorname{H}^{p}(\c/I,0;(M^{*}\otimes M)^{I}) \cong \operatorname{H}^{p}(\c,\c_{\0};M^{*}\otimes M)
\] for all $p \geq 0.$  Then by Theorem~\ref{T:finitegen2}, $\operatorname{H}^{\bullet}(\c, \c_{\0};M^{*}\otimes M)$ is finitely generated as a module for
\[
\operatorname{H}^{\bullet}(\c/I, 0;\C) \cong S((\c/I)_{\1})^{*})\cong S((\c_{\1}/\mathbb{C}x)^{*}).
\]

Let 
\begin{align*}
N &=\operatorname{H}^{\bullet}(\c, \c_{\0};M^{*}\otimes M) \otimes_{S(\c_{\1}^{*})} S(\c_{\1}^{*})/\left( (\c_{\1}/\mathbb{C}x)^{*}S(\c_{\1}^{*})\right)\\
&\cong \operatorname{H}^{\bullet}(\c, \c_{\0};M^{*}\otimes M) \otimes_{S((\c_{\1}/\mathbb{C}x)^{*})} \mathbb{C}.
\end{align*}  From commutative algebra (e.g.\ \cite[Chapter 3, Ex.\ 19(vi)]{atiyah}) we have that 
\[
\operatorname{Supp}(N):=\left\{ \mathfrak{m} \in 
{\mathcal V}_{(\c,\c_{\0})}(\C)\:\vert\:  N_{\mathfrak{m}} \neq 0 \right\}={\mathcal V}_{(\c,\c_{\0})}(M) \cap \res^{*}(
{\mathcal V}_{(\langle x \rangle,\langle x \rangle_{\0})}(\mathbb{C})).
\]
However, since $\operatorname{H}^{\bullet}(\c, \c_{\0 };M^{*}\otimes M)$ is a finitely 
generated $S((\c_{\1}/\mathbb{C}x)^{*})$-module, $N$ is finite dimensional.  This implies $\operatorname{Supp}(N)$ is a finite set.  On the other hand, 
$\operatorname{Supp}(N)$ is conical.  Therefore, $\operatorname{Supp}(N)=\{0 \}.$  That is, since $ \res^{*}({\mathcal V}_{(\langle x \rangle,\langle x \rangle_{\0})}(\mathbb{C}))=\C x,$ $x \not\in {\mathcal V}_{(\c,\c_{\0})}(M).$
\end{proof}

\subsection{}\Label{SS:complexity} \setcounter{Df}{0} One can use the rank variety description of 
${\mathcal V}_{(\c,\c_{\0})}(M)$ along with Lemmas~\ref{L:dualizingprojs} and~\ref{L:irrepsinFL} to prove the analogue of well known results in the theory of support varieties.  In most cases the classical arguments from finite groups and restricted Lie algebras apply here as well.  For these standard arguments we refer the reader to \cite{benson, FP}. 

\begin{corollary}\Label{C:rankcor1} Let $M$ be a finite dimensional $\c$-supermodule. If $\dim M_{\0} \neq \dim M_{\1}$ (e.g.\ when $M$ is odd dimensional), then 
\[
{\mathcal V}_{(\c,\c_{\0})}(M)={\mathcal V}_{(\c,\c_{\0})}(\mathbb{C}).
\]
\end{corollary}

\begin{proof} This follows from the rank variety description of 
${\mathcal V}_{(\c,\c_{\0})}(M).$   Namely, fix $0 \neq x \in \c_{\1}.$   By Proposition~\ref{L:irrepsofQ(1)}, the projective indecomposable $\langle x \rangle$-supermodules $P$ satisfy $\dim P_{\0}=\dim P_{\1} = 1.$   Thus if $M$ is a projective $\langle x  \rangle$-supermodule, then  $\dim M_{\0} = \dim M_{\1}.$  As this contradicts what is assumed of $M,$ $x \in \mathcal{V}_{\c}^{\operatorname{rank}}(M).$  Therefore ${\mathcal V}_{(\c,\c_{\0})}(M)=\mathcal{V}_{\c}^{\operatorname{rank}}(M)=\c_{\1}={\mathcal V}_{(\c,\c_{\0})}(\mathbb{C}).$  
\end{proof}

Before continuing, let us recall the following basic notions in homological algebra.  If $\mathfrak{a}$ is a 
Lie superalgebra, then the \emph{complexity} of a finite dimensional $\mathfrak{a}$-supermodule $M$ is the rate of growth of a minimal projective resolution of $M$ (cf.\ \cite[Sec. 5.3]{benson}). Given a finite dimensional $\mathfrak{a}$-supermodule $M,$ let $P(M)$ denote the projective cover of $M.$  Recalling that we assume $P(M)$ maps onto $M$ by an even supermodule homomorphism,  one defines $\Omega(M)=\Omega^{1}(M)$ to be the kernel of the aforementioned homomorphism.  For $n> 1$ recursively define $\Omega^{n}(M)$ to be $\Omega(\Omega^{n-1}(M)).$  

The following theorem taking together with Proposition~\ref{L:tensorproduct} shows that support varieties for finite dimensional $\c$-supermodules satisfy the desirable properties of a support variety theory.

\begin{theorem}\label{C:tensorproducttheorem}  If $M,$ $M_{1},$ $M_{2},$ $M_{3},$ and $N$ are finite dimensional $\c$-supermodules, then 
\begin{itemize}
\item[(a)] The $\text{complexity of $M$}=\dim \mathcal{V}_{(\c,\c_{\0})}(M).$
\item[(b)] $M$ is projective if and only if $\mathcal{V}_{(\c,\c_{\0})}(M)=\{0 \}$. 
\item [(c)] ${\mathcal V}_{(\c,\c_{\0})}(M_{1}\oplus M_{2}) = {\mathcal V}_{(\c,\c_{\0})}(M_{1}) \cup{\mathcal V}_{(\c,\c_{\0})}(M_{2}).$
\item [(d)] If 
\[
0 \to M_{1} \to M_{2} \to M_{3} \to 0
\]
is a short exact sequence of $\c$-supermodules, then 
\[
{\mathcal V}_{(\c,\c_{\0})}(M_{i}) \subseteq {\mathcal V}_{(\c,\c_{\0})}(M_{j}) \cup{\mathcal V}_{(\c,\c_{\0})}(M_{k}), 
\]
where $\{i,j,k \}=\{1,2,3 \}.$
\item [(e)]  For any $n \geq 1,$  $\mathcal{V}_{(\c,\c_{\0})}(\Omega^{n}(M))=\mathcal{V}_{(\c,\c_{\0})}(M).$
\item [(f)] If ${\mathcal V}_{(\c,\c_{\0})}(M_{1}) \cap{\mathcal V}_{(\c,\c_{\0})}(M_{2})=\{0 \},$ then for any finite dimensional $\c$-supermodule $N$ the short exact sequence 
\[
0 \to M_{1} \to N \to M_{2} \to 0
\]
splits.
\end{itemize}
\end{theorem}

\begin{proof}   Parts (a) and (b) are proved by the same arguments used for restricted Lie algebras \cite[Propostions 1.5, 3.2]{FP}. Part (c) follows immediately from the rank variety description and part (d) follows from the rank variety description along with Lemma~\ref{L:irrepsofQ(1)}.  Part (e) follows from part (d) and the definition of $\Omega^{n}(M).$  Part (f) follows just as in the classical case using parts (a) and (d).
\end{proof}

Just as for finite groups, one has $\operatorname{H}^{p}(\c,\c_{\0};\C)\cong \Hom_{\c}(\Omega^{p}(\C),\C).$  Given $\zeta \in \operatorname{H}^{p}(\c,\c_{\0 };\C)$ we let 
\begin{equation}\label{E:Lzetadef}
L_{\zeta}=\Ker (\zeta: \Omega^{p}(\C) \to \C) \subseteq \Omega^{p}(\C).
\end{equation}  As in the classical theory of support varieties (cf.\ \cite[Lemma 4.2]{FP}), the importance of the supermodule $L_{\zeta}$ is that one can explicitly calculate its support variety.  In particular, these supermodules allow one to prove the following realization theorem.

\begin{theorem}\label{T:realization}  Let $X$ be a closed conical subvariety of $\mathcal{V}_{(\c,\c_{\0})}(\C).$  There is a finite dimensional $\c$-supermodule $M$ such that 
\[
\mathcal{V}_{(\c,\c_{\0})}(M) = X.
\]
\end{theorem}

\begin{proof}  First, given $\zeta \in \operatorname{H}^{p}(\c,\c_{\0};\C),$ one can compute the support variety of $L_{\zeta}$ as follows.  From the graded version of Schanuel's Lemma, $\Omega^{p}(\C) \cong P \oplus \Omega_{\langle x \rangle}^{p}(\C)$ as $\langle x \rangle$-supermodules, where $\Omega_{\langle x \rangle}^p(\C)$ denotes $\Omega^{p}(\C)$ for the trivial $\langle x \rangle$-supermodule and  $P$ is some projective $\langle x \rangle$-supermodule.  By Lemma~\ref{L:irrepsofQ(1)}, $\Omega_{\langle x \rangle}^{p}(\C) \cong \C.$  With these facts in hand one can argue just as for finite groups and restricted Lie algebras (cf.\ \cite[Lemma 4.1]{FP}) to prove 
\begin{equation}\label{E:Lzetavariety}
\mathcal{V}_{(\c,\c_{\0})}(L_{\zeta}) = \mathcal{Z}(\zeta):=\left\{ x \in \c_{\1} \:\vert\:  \zeta (x) =0 \right\}.
\end{equation}

Now let $I =(\zeta_{1}, \dotsc , \zeta_{t})\subseteq \operatorname{H}^{\bullet}(\c,\c_{\0};\C) $ be a homogeneous ideal which defines $X.$  That is, $\zeta_{1}, \dotsc , \zeta_{t}$ are homogeneous elements of $\operatorname{H}^{\bullet}(\c,\c_{\0};\C)$ and
\[
X=\mathcal{Z}(\zeta_{1}) \cap \dotsb \cap \mathcal{Z}(\zeta_{t}).
\]
Applying \eqref{E:Lzetavariety} and Proposition~\ref{L:tensorproduct}(a) one has 
\begin{align*}
\mathcal{V}_{(\c,\c_{\0})}(L_{\zeta_{1}} \otimes \dotsb \otimes L_{\zeta_{t}})= \mathcal{Z}(\zeta_{1}) \cap \dotsb \cap \mathcal{Z}(\zeta_{t})= X.
\end{align*}
\end{proof}

We also have the following connectedness result.

\begin{theorem}\label{T:connectedness} Suppose that $M$ is a finite dimensional $\c$-supermodule such that 
\[
\mathcal{V}_{(\c,\c_{\0})}(M)=X_{1} \cup X_{2},
\]
where $X_{1}$ and $X_{2}$ are nonzero closed conical subvarieties with $X_{1}\cap X_{2}=\{0 \}.$  Then there are finite dimensional $\c$-supermodules $M_{1}$ and $M_{2}$ such that $\mathcal{V}_{(\c,\c_{\0})}(M_{1})=X_{1},$ $\mathcal{V}_{(\c,\c_{\0})}(M_{2})=X_{2},$ and $M \cong M_{1}\oplus M_{2}.$

In particular, if $M$ is indecomposable, then the projectivization of the conical variety $\mathcal{V}_{(\c,\c_{\0})}(M)$ is connected.
\end{theorem}

\begin{proof}  One can argue just as in the case of finite groups (cf.\ \cite{carlson}) using \eqref{E:Lzetavariety}, Proposition~\ref{L:tensorproduct}, and Theorem~\ref{C:tensorproducttheorem}.
\end{proof}
\section{Defect, Atypicality and Superdimension} 

\subsection{Defect}\Label{SS:defect} \setcounter{Df}{0}Let ${\mathfrak g}$ be a classical Lie superalgebra and ${\mathfrak t}$ be 
a maximal torus contained in $\g_{\0}$. Let $\Phi$ be the set of roots with respect 
to ${\mathfrak t}$. We have that $\Phi=\Phi_{\0}\cup \Phi_{\1}$ 
where $\Phi_{\0}$ (resp.\ $\Phi_{\1}$) is the set of even roots (resp.\ odd roots). 
The positive roots will be denoted by $\Phi^{+}$ and the negative roots by $\Phi^{-}$. 
Set $\Phi_{\0}^{\pm}=\Phi_{\0}\cap \Phi^{\pm}$ and $\Phi_{\1}^{\pm}=\Phi_{\1}\cap \Phi^{\pm}$. 
If $\g$ is a basic classical Lie superalgebra (cf.\ Section~\ref{S:prelims}), then Kac and Wakimoto \cite[\S2]{kacwakimoto} 
define the {\em defect of $\g$}, denoted by $\text{def}({\mathfrak g}),$ to be the 
dimension of a maximal isotropic subspace in the ${\mathbb R}$-span of $\Phi$.

The defects for the various simple basic classical Lie superalgebras are given as follows 
\cite[Section 4]{dufloserganova}: $\text{def}(\mathfrak{sl}(m|n))=\text{min}(m,n)$,  $\text{def}(\mathfrak{psl}(n|n))=n,$
$\text{def}(\mathfrak{osp}(2m+1|2n))=\text{def}(\mathfrak{osp}(2m|2n))=\text{min}(m,n)$, and the exceptional
Lie superalgebras $D(2,1;\alpha),$ $G(3),$ and $F(4)$ all have defect 1.  The 
following theorem demonstrates that one can realize this numerical 
defect using relative cohomology, support varieties, or the detecting subalgebra $\c.$

\begin{theorem}\label{T:defect} Let ${\mathfrak g}$ be a basic classical simple Lie superalgebra. Moreover, assume that 
${\mathfrak g}\ncong \mathfrak{psl}(n|n)$. The following 
numbers are equal. 
\begin{itemize}
\item[(a)] $\operatorname{def}({\mathfrak g})$;
\item[(b)] $\dim \operatorname{H}^{\bullet}(\g,\g_{\0};\C)$;  
\item[(c)] $\dim {\mathcal V}_{(\g,\g_{\0})}({\mathbb C})$.
\end{itemize}
Furthermore, if $\g$ is stable and polar, then the above are also equal to the following numbers.
\begin{itemize}
\item[(d)] $\dim \operatorname{H}^{\bullet}({\mathfrak a},{\mathfrak a}_{\0};\C),$ where $\mathfrak{a}=\e$ or $\c$; 
\item[(e)] $\dim {\mathcal V}_{(\mathfrak{a},\mathfrak{a}_{\0})}({\mathbb C}),$ where $\mathfrak{a}=\e$ or $\c$;
\item[(f)] $\dim {\mathcal V}^{\operatorname{rank}}_{\c}(\C);$ 
\item[(g)] $\dim \c_{\1}.$  
\end{itemize} 
\end{theorem} 

\begin{proof}  The equality of (a)--(c) follows by verifying that the defects listed above equal the dimensions presented in Table~1 of the Appendix.   That these coincide with the dimensions listed in (d)--(g) follow from the results 
in Sections~\ref{S:alggeo} and~\ref{S:morecohom}. 
\end{proof} 

The above theorem indicates that one could naturally extend the definition of defect 
cohomologically to arbitrary classical Lie superalgebras by setting $\text{def}({\mathfrak g})$ to 
be the Krull dimension of $\operatorname{H}^{\bullet}(\g,\g_{\0};\C)$.  In particular, the new definition would allow one to define 
the notion of defect for the simple classical Lie superalgebras of type $P(n)$ and $Q(n).$ 
Also, let us remark that the definition given by Kac and Wakimoto would differ 
from this new definition for the Lie superalgebra $\mathfrak{psl}(n|n)$. Namely, $\operatorname{def}(\mathfrak{psl}(n|n))=\operatorname{def}(\mathfrak{gl}(n|n))=n,$ whereas  
\[
\dim \operatorname{H}^{\bullet}(\mathfrak{psl}(n|n),\mathfrak{psl}(n|n)_{\0};\C)=n+1 \neq \dim \operatorname{H}^{\bullet}(\mathfrak{gl}(n|n),\mathfrak{gl}(n|n)_{\0};\C)=n.
\]

\subsection{Atypicality}\Label{SS:atypicality} \setcounter{Df}{0} Let $\g$ be a basic classical Lie superalgebra as above and let $(-,-)$ denote the bilinear form.  Let $\lambda\in {\mathfrak t}^{*}$ be a weight. The \emph{atypicality} of $\lambda$ is the maximal number of linearily independent, mutually orthogonal, positive isotropic roots $\alpha \in  \Phi^{+}$ such that $ (\lambda+\rho,\alpha)=0,$
%$$\operatorname{atyp}(\lambda)=|\{ \alpha\in \Phi_{\1}^{+} \:\vert\:\ (\lambda+\rho,\alpha)=0\}|$$ 
where $\rho=\frac{1}{2}(\sum_{\alpha\in \Phi_{\0}^{+}} \alpha-
\sum_{\alpha\in \Phi_{\1}^{+}} \alpha)$. Note that $\text{atyp}(\lambda)\leq \text{def}({\mathfrak g})$. 
Let $L(\lambda)$ be a simple finite-dimensional ${\mathfrak g}$-supermodule with highest weight $\lambda$.  
The atypicality of $L(\lambda)$, denoted by $\text{atyp}(L(\lambda))$, is defined to be $\text{atyp}(\lambda)$. 
We present the following strong conjecture. 

\begin{conj}\label{C:atypconjecture} Let $\g$ be a simple basic classical Lie superalgebra which is stable and polar and let $L(\lambda)$ be a finite dimensional simple $\g$-supermodule. 
Then $$\operatorname{atyp}(L(\lambda))=\dim {\mathcal V}^{\operatorname{rank}}_{{\mathfrak e}}(L(\lambda)).$$ 
\end{conj} 

There is the following evidence in favor of the conjecture.  Using results in \cite{kacwakimoto} one can verify the conjecture for all $\g$ which have defect one.
The conjecture also holds in the case when the atypicality of the simple $\g$-supermodule $L(\lambda)$ is zero. Namely, $\text{atyp}(L(\lambda))=0$ implies $L(\lambda)$ is a projective $U({\mathfrak g})$-supermodule (cf.\ \cite[Theorem 1]{Kacnote}), thus projective 
as a $U(\c)$-supermodule. This in turn implies ${\mathcal V}^{\operatorname{rank}}_{{\mathfrak e}}
(L(\lambda))=\{0\}$.  Finally, in the case of $\g=\gl,$ one can use translation functor techniques from \cite{dufloserganova} to prove $\operatorname{atyp}(L(\lambda))\leq \dim {\mathcal V}^{\operatorname{rank}}_{{\mathfrak e}}(L(\lambda))$ for any finite dimensional simple $\g$-supermodule $L(\lambda).$ 

\subsection{Kac-Wakimoto Conjecture} \setcounter{Df}{0}Let $M$ be a supermodule. The {\em superdimension} 
of $M$ is defined to be $\text{sdim }M=\dim M_{\0}-\dim M_{\1}$. Kac and Wakimoto 
give the following conjecture relating the superdimension of simple finite-dimensional ${\mathfrak g}$-supermodules 
with the atypicality of the module and the defect of ${\mathfrak g}$. 

\begin{conj} \cite[Conjecture 3.1]{kacwakimoto} Let $\g$ be a simple basic classical Lie superalgebra 
and $L(\lambda)$ be a finite-dimensional simple ${\mathfrak g}$-module. Then 
$\operatorname{sdim }L(\lambda)=0$ if and only if 
$\operatorname{atyp}(L(\lambda))< \operatorname{def}({\mathfrak g})$.  
\end{conj} 

We will now indicate how our conjecture given in Section~\ref{SS:atypicality} and the 
Kac-Wakimoto Conjecture are interrelated. On the one hand, the validity of 
our conjecture implies one direction of the Kac-Wakimoto Conjecture. 
Suppose that $\operatorname{sdim }L(\lambda)\neq 0.$  Then 
by Corollary~\ref{C:rankcor1}, 
$\dim {\mathcal V}^{\operatorname{rank}}_{{\mathfrak e}}(L(\lambda))=\dim {\mathfrak e}_{\1}$. 
So by using Conjecture~\ref{C:atypconjecture} and Theorem~\ref{T:defect}, 
$\operatorname{atyp}(L(\lambda))=\text{def}({\mathfrak g})$.  Note that this implication of the Kac-Wakimoto conjecture was recently verified for the simple classical contragradient Lie superalgebras \cite[Lemma 7.1]{dufloserganova}.
On the other hand, assume that the Kac-Wakimoto Conjecture is valid. Then our conjecture 
is true for all simple ${\mathfrak g}$-modules with atypicality equal to $\text{def}({\mathfrak g})$. 
If $\text{atyp}(L(\lambda))=\text{def}({\mathfrak g})$ then $\operatorname{sdim}L(\lambda)\neq 0$. 
Therefore, by Corollary~\ref{C:rankcor1}, ${\mathcal V}^{\operatorname{rank}}_{{\mathfrak e}}(L(\lambda))
={\mathfrak e}_{\1}$.

\section{Appendix:  Superalgebra Data and Sample Computations}\Label{S:cohomorings}  In this appendix we record various data for the simple Lie superalgebras of classical type and other related Lie superalgebras.  We begin by quickly reviewing the structure of the simple classical Lie superalgebras; for more details on their definition and structure, including a precise description of $\fg_{\1}$ as a $\fg_{\0}$-supermodule, we refer the reader to \cite[Chap. 2]{K}.  We next tabulate their cohomology rings.  We then define the generic semisimple elements $x_0$ that one can use in both the Dadok-Kac and Luna-Richardson reduction arguments, and provide the data necessary to determine when $x_0$ is regular. We also provide explicit descriptions of the detecting subalgebras $\e$ and $\c.$  Finally we give examples of both of the Dadok-Kac and Luna-Richardson computations in type $A$. 

\subsection{The Lie superalgebras of type A}\Label{SS:typeA}   Let $\fg=\mathfrak{gl}(m|n),$ the Lie superalgebra of $(m+n) \times (m+n)$ matrices over $\C$ with $\Z_2$-grading given by setting the degree of the matrix unit $E_{i,j}$ to be $\0$ if $1 \leq i,j \leq m$ or $m+1 \leq i,j \leq  m+n$ and $\1$ otherwise.  The bracket is the super-commutator bracket, 
\[
[a,b] =ab-(-1)^{\p{a}\p{b}}ba,
\] where $a,b \in \fg$ are assumed to be homogenous.  The bracket of general elements is obtained by bilinearity.  It is straightforward to verify that 
\[
\fg_{\0 } \cong \mathfrak{gl}(m) \oplus \mathfrak{gl}(n) \text{ and } \fg_{\1} \cong V_{m} \boxtimes V_{n}^{*} \oplus V_{m}^{*} \boxtimes V_{n},
\] where the first isomorphism is as Lie algebras, the second isomorphism is as $\fg_{\0}$-modules, and where $V_{m}$ (resp.\ $V_{n}$) denotes the natural $\mathfrak{gl}(m)$-module (resp.\ $\mathfrak{gl}(n)$-module).  For convenience, we shall say that $\gl$ is a Lie superalgebra of type $\widehat A(m-1,n-1)$.

Let $m,n \geq 1$, $m \neq n$ and let $\fg=\slmn$ be the simple Lie superalgebra of type $A(m-1,n-1).$  That is, $\fg$ is the Lie subsuperalgebra of $\gl$ of all matrices with supertrace equal to zero.  Then 
\[
\fg_{\0 } \cong \mathfrak{sl}(m) \oplus \mathfrak{sl}(n) \oplus \C \text{ and } \fg_{\1} \cong V_{m} \boxtimes V_{n}^{*}\boxtimes \C \oplus V_{m}^{*} \boxtimes V_{n} \boxtimes \C,
\] where $\C$ denotes the one-dimensional trivial Lie algebra and the trivial representation, respectively.

Let $m=n \geq 2$ and let $\fg=\slnn$ be the simple Lie superalgebra of type $A(n-1,n-1).$  That is, if $\widetilde{\fg}=\mathfrak{sl}(n|n)$ is the Lie subsuperalgebra of $\gl$ of all matrices with supertrace equal to zero, then $\fg$ is the quotient of $\widetilde{\fg}$ by the one-dimensional ideal spanned by the identity matrix.  Then one obtains 
\[
\fg_{\0 } \cong \mathfrak{sl}(n) \oplus \mathfrak{sl}(n) \text{ and } \fg_{\1} \cong V_{n} \boxtimes V_{n}^{*} \oplus V_{n}^{*} \boxtimes V_{n}.
\]

\subsection{The simple Lie superalgebra of type P}\Label{SS:typeP} Let $\fg$ be the simple Lie superalgebra of type $P(n-1)$ where $n\ge 3$, consisting of $2n \times 2n$ matrices of the form
\begin{equation}\Label{E:Pmatrix}\left( 
\begin{array}{c|c}
A&B\\\hline
C&-A^t
\end{array}
\right),
\end{equation}
where $A, B$ and $C$ are $n \times n$ matrices over $\C$ with $\tr A=0$, $B$ symmetric, and $C$ skew-symmetric.  Then
\[
\fg_{\0 } \cong \mathfrak{sl}(n) \text{ and } \fg_{\1} \cong S^2 V_n \oplus \wedge^{2} V_n^*,
\]
where $V_n$ is the natural $\mathfrak{sl}(n)$-module.

\subsection{The Lie superalgebras of type Q}\Label{SS:typeQ} Let $\fg=\mathfrak{q}(n)$ denote the Lie subsuperalgebra of $\mathfrak{gl}(n|n)$ consisting of $2n \times 2n$ matrices of the form   
\begin{equation}\Label{E:Qmatrix}\left( 
\begin{array}{c|c}
A&B\\\hline
B&A
\end{array}
\right),
\end{equation}
where $A$ and $B$ are both $n \times n$ matrices over $\C$.  It is again straightforward to verify that 
\[
\fg_{\0 } \cong \mathfrak{gl}(n) \text{ and } \fg_{\1} \cong \mathfrak{gl}(n),
\] where we mean $\fg_{\1}$ is the adjoint representation for $\fg_{\0}.$ For convenience we shall say that $\fq(n)$ is a Lie superalgebra of type $\widehat Q(n-1)$.

Let $\fg$ denote the simple Lie superalgebra of type $Q(n-1)$ where $n \geq 3.$  That is, $\fg$ is the Lie subsuperalgebra of $\mathfrak{q}(n)$ of all matrices of the form \eqref{E:Qmatrix} such that $A$ and $B$ have trace zero. Then clearly 
\[
\fg_{\0 } \cong \mathfrak{sl}(n) \text{ and } \fg_{\1} \cong \mathfrak{sl}(n).
\] 

\subsection{The Lie superalgebras of type BCD}\Label{SS:typeBCD}  Let $\fg=\mathfrak{osp}(m|n).$  That is, let $V$ denote a superspace with $\dim V_{\0}=m$ and $\dim V_{\1}=n$ and equipped with a fixed nondegenerate supersymmetric bilinear form $\beta$ such that with respect to this form $V_{\0}$ and $V_{\1 }$ are orthogonal subspaces.  In particular, note that the form is skew symmetric on $V_{\1 }$ so necessarily $n$ is even.  One defines $\g$ to be the Lie subsuperalgebra of $\gl$ given by
\[
\mathfrak{osp}(m|n)= \setof{x \in \mathfrak{gl}(V)}{\beta(x(v),w)=-(-1)^{\p{x}\;\p{v}}\beta(v,x(w))}
\]
Then 
\[
\fg_{\0 } \cong \mathfrak{so}(m) \oplus \mathfrak{sp}(n) \text{ and } \fg_{\1} \cong  V_{m} \boxtimes V_{n},
\] where $V_{m}$ (resp.\ $V_{n}$) is the natural representation for $\mathfrak{so}(m)$ (resp.\ $\mathfrak{sp}(n)$).  

The classification types of the simple Lie algebras $\osp(m|n)$ are as follows: type $B(m,n)=\osp(2m+1|2n)$, $m\ge 0, n>0$; type $D(m,n)=\osp(2m|2n)$, $m\ge 2, n>0$; type $C(n)=\osp(2|2n-2)$, $n\ge 2$.

\subsection{The simple Lie superalgebra of type $D(2,1;\alpha)$}\Label{SS:typeD21alpha} Let $\fg$ denote the simple Lie subsuperalgebra of type $D(2,1;\alpha)$ ($\alpha \in \C$) as described in \cite{K}.  Then 
\[
\fg_{\0 } \cong \mathfrak{sl}(2) \oplus \mathfrak{sl}(2) \oplus \mathfrak{sl}(2) \text{ and } \fg_{\1} \cong V_{2} \boxtimes V_{2} \boxtimes V_{2},
\] where $V_{2}$ denotes the natural $\mathfrak{sl}(2)$-module.

\subsection{The simple Lie superalgebra of type $G(3)$}\Label{SS:typeG3}  Let $\fg$ denote the simple Lie superalgebra of type $G(3).$  Then 
\[
\fg_{\0 } \cong \mathfrak{a} \oplus \mathfrak{sl}(2) \text{ and } \fg_{\1} \cong W_{7} \boxtimes V_{2},
\] where $\mathfrak{a}$ denotes the simple Lie superalgebra of type $G_{2},$ $W_{7}$ denotes the nontrivial $\mathfrak{a}$-module of minimal dimension and $V_{2}$ denotes the natural $\mathfrak{sl}(2)$-module.

\subsection{The simple Lie superalgebra of type $F(4)$}\Label{SS:typeF4}  Let $\fg$ denote the simple Lie superalgebra of type $F(4).$  Then 
\[
\fg_{\0 } \cong \mathfrak{a} \oplus \mathfrak{sl}(2) \text{ and } \fg_{\1} \cong \operatorname{spin}_{7} \boxtimes V_{2},
\] where $\mathfrak{a}$ denotes the simple Lie superalgebra of type $B_{3},$ $\operatorname{spin}_{7}$ denotes the $\mathfrak{a}$-module $\operatorname{spin}_{7}$ and $V_{2}$ denotes the natural $\mathfrak{sl}(2)$-module.

\subsection{Cohomology rings}\Label{SS:cohomrings} One of our main objects of interest is the relative cohomology ring $\operatorname{H}^\bullet(\fg,\fg_{\0};\C)$, where $\fg$ is one of the classical Lie superalgebras described in the previous subsections. According to Theorem~\ref{T:cohomring}, $\operatorname{H}^\bullet(\fg,\fg_{\0};\C)\cong S(\fg_{\1}^*)^{G_{\0}}$, and so we can exhibit these rings simply by referring to known calculations in invariant theory. They turn out to be polynomial rings in every case.

In Table~\ref{T:cohomrings} we list the Krull dimension of each relative cohomology ring, the degrees of its generators in the $\Z$-grading (recall that the $\Z_{2}$-grading is obtained by reducing the $\Z$-grading modulo $2$), and a reference where the result can be found. Given the rich history of invariant theory, we make no claim to our references as the original source of these calculations. By convention, set $r:=\min(m,n)$.

\renewcommand{\arraystretch}{1.2}
\begin{table}[htdp]
\caption{Cohomology rings for classical Lie superalgebras}
%\moveleft 20pt
\vbox{
\begin{center}
\begin{tabular}{ccccccl}
$\fg$ &  $\dim S(\fg_{\1}^*)^{G_{\0}}$ & Degrees of Generators & Reference \\ \hline
$\gl$ & $r$ & $2,4,\dots,2r$ & \cite[Sec.\ 2.1.2]{fuks} \\
$\slmn$, $m\ne n$& $r$ & $2,4,\dots,2r$ & \cite{Gru:97} \\
%\raisebox{.5ex}{$m\ne n$} \\
$\slnn$ &  $n+1$ & $2,4,\dots,2n-2, n, n$ & Section \ref{SS:LR}\\
$\mathfrak{osp}(2m+1|2n)$ & $r$ & $4,8,\dots,4r$ & \cite[Tables II, III]{Kac2}\\
$\mathfrak{osp}(2m|2n)$, $m>n$& $n$ & $4,8,\dots, 4n$ & \cite[Tables II, III]{Kac2}\\
%\raisebox{.5ex}{$m>n$} \\
$\mathfrak{osp}(2m|2n)$, $m\le n$& $m$ & $4,8,\dots,4m-4,2m$ &  \cite[Tables II, III]{Kac2}\footnotemark \\
%\raisebox{.5ex}{$m\le n$} \\
$P(n-1)$, $n=2l+1$& $l+1$ & $4,8,\dots,4l,n$ & \cite{Gru:00, Sch:78}\\
%\raisebox{.5ex}{$n=2l+1$} \\
$P(n-1)$, $n=2l$& $l+1$ & $4,8,\dots,4l-4,l,n$ & \cite{Gru:00, Sch:78}\\
%\raisebox{.5ex}{$n=2l$} \\
$\widehat Q(n-1)$ & $n$ &  $1, 2,\dots,n$ & \cite{Weyl:39} \\
$Q(n-1)$ & $n-1$ & $2, 3,\dots,n$ & \cite{Weyl:39} \\
$D(2,1;\a)$ & 1 & $4$ & \cite{MeyerWallach}\\
$G(3)$ & 1 & $4$ &\cite[Table IV]{Kac2}\\
$F(4)$ & 1 &  $4$ &\cite[Table III]{Kac2}
\end{tabular}
\end{center}
}
\Label{T:cohomrings}
\end{table}%

\subsection{Explicit detecting subalgebras} \Label{SS:explicitdetecting}
In this and the next two subsections, we give explicit descriptions of the detecting subalgebras (cf.\ Section \ref{SS:detecting}) for each classical Lie superalgebra $\fg$. Recall that we begin by defining a generic element $x_0\in \fg_{\1}$; this will be a certain linear combination of sums of positive and negative odd root vectors. 
Using the following theorem of Dadok and Kac \cite[Proposition 1.2]{dadokkac}, we will deduce that the $G_{\0}$-orbit of $x_0$ is closed; that is, $x_0$ is semisimple.

\begin{theorem}\Label{T:closedorbits}  Suppose $v_{\mu_{1}}, \dotsc , v_{\mu_{s}}$ is a set of nonzero weight vectors in a $G$-module $V$ such that 

\begin{enumerate}
\item the weights $\mu_{1}, \dotsc ,\mu_{s}$ are distinct;
\item 0 is an interior point of the convex hull of $\mu_{1}, \dotsc ,\mu_{s};$
\item $\mu_{i}-\mu_{j}$ is not a root if $i \neq j.$
\end{enumerate}  Then the $G$ orbit of $v_{\mu_{1}}+\dotsb +v_{\mu_{s}}$ is closed in $V.$
\end{theorem}

Furthermore, in all but two cases, the subspace $\fc_{x_0}$ is a Cartan subspace of $\fg_{\1}$, so the action of $G_{\0}$ on $\fg_{\1}$ is polar (cf.\ Section~\ref{SS:dadokkac}). One can then use the Dadok-Kac theorem \cite{dadokkac} to compute the cohomology ring $S(\fg_{\1}^*)^{G_{\0}}$ as in Section~\ref{SS:dadokkac}. An example of these computations is given in Section~\ref{SS:DK}.

\addtocounter{footnote}{0}\footnotetext{In \cite[Table II]{Kac2}, for $SL_n\otimes SO_m$, $1\le n < m$, $d$ should be $2n$, not $2m$, and for $SL_n\otimes Sp_m$, $2\le n\le m$, $n$ even, $d$ should be $n$, not $m$. And in \cite[Table III]{Kac2}, for $SO_n\otimes Sp_m$, $2<n\le m$, $n$ even, $d_1$ should be $n$, not $m$.}

In the remaining two cases, the action is stable, and we can apply the Luna-Richardson theorem \cite{luna} to compute $S(\fg_{\1}^*)^{G_{\0}}$ as in Section~\ref{SS:stableactions}.  An example computation is carried out in Section~\ref{SS:LR}.

We begin with the definitions.  Let $\fg$ be one of the classical Lie superalgebras discussed in Sections~\ref{SS:typeA}--\ref{SS:typeF4} other than types $\widehat{Q}$ and $Q$. Let $\Omega$ be the set of odd positive roots defined in Table~\ref{T:Omega}. (As before, we set $r=\min(m,n)$. Notation for the roots follows \cite{K}.) Fix an element $x_0=\sum_{\a\in\Omega} d_\a(x_\a+x_{-\a}) \in \fg_{\1}$, where $(d_\a)_{\a\in\Omega}$ is a vector of complex numbers which is not in the zero locus of the Jacobian $J$ of Table \ref{T:polarstable}. (For $\slnn$ and $P(n-1)$ we assume the $d_\a$ are nonzero with distinct squares.) The set $\Omega \cup -\Omega$ satisfies the conditions of Theorem \ref{T:closedorbits} and so $x_0$ is semisimple.

\renewcommand{\arraystretch}{1.2}
\begin{table}[htdp]
\caption{Sets $\Omega$ defining $x_0,\ \fc_{\1}$, and $\fe_{\1}$ for classical Lie superalgebras}
\begin{center}
\begin{tabular}{cc}
$\fg$ & $\Omega$  \\ \hline
$\gl$ & $\setof{\ep_i-\ep_{m+i}}{1\le i\le r}$ \\
$\slmn,\ m\ne n$ & $\setof{\ep_i-\ep_{m+i}}{1\le i\le r}$ \\
$\slnn$ & $\setof{\ep_i-\ep_{n+i}}{1\le i\le n}$ \\
$\mathfrak{osp}(2m+1|2n)$ & $\setof{\ep_i-\d_i}{1\le i\le r}$\\
$\mathfrak{osp}(2m|2n)$ & $\setof{\ep_i-\d_i}{1\le i\le r}$\\
$P(n-1),\ n$ even & $\setof{\ep_i+\ep_{n+1-i}}{1\le i\le \frac{n}{2}}$\\
$P(n-1),\ n$ odd & $\setof{\ep_i+\ep_{n+1-i}}{1\le i\le \frac{n+1}{2}}$\footnotemark\\
%$\widehat Q(n-1)$ & $\{0,\dots,0\}$\footnotemark \\
%$Q(n-1)$ & $\{0,\dots,0\}$\addtocounter{footnote}{-1}\footnotemark \\
$D(2,1;\a)$ & $\{ \ep_1+\ep_2+\ep_3\}$\\
$G(3)$ & $\{ \ep_1+\d \}$\\
$F(4)$ & $\{ \frac12(\ep_1+\ep_2+\ep_3+\ep_4) \}$
\end{tabular}
\end{center}
\Label{T:Omega}
\end{table}%

Define subspaces 
\begin{align}
\fc_{\1} &= \setof{\sum_{\a\in\Omega} c_\a(x_\a+x_{-\a})}{c_\a\in\C}, \Label{E:cdef}\\
\fe_{\1} &= \setof{\sum_{\a\in\Omega} (u_\a x_\a+v_\a x_{-\a})}{u_\a, v_\a\in\C}  \Label{E:edef}
\end{align}
of $\fg_{\1}$. Evidently $x_0\in\fc_{\1}\subset\fe_{\1}$.

Now assume we are in one of the cases $\widehat Q(n-1)$ or $Q(n-1)$. Recall from Section \ref{SS:typeQ} that an element $x\in\fg_1$ is determined by an $n\times n$ matrix $B$ as in \eqref{E:Qmatrix}. When $B$ is diagonal (having trace zero in type $Q(n-1)$), $x$ lies in the 0-root space. We take $x_0$ to be an element of $\fg_{\1}$ whose $B$ part is a diagonal matrix with $n$ distinct diagonal entries.  From the classical theory of $\mathfrak{gl}(n)$ and $\mathfrak{sl}(n)$ one knows that $x_0$ is generic.  We let $\fc_{\1}$ be the set of all $x\in\fg_{\1}$ whose $B$ matrix is diagonal (and trace zero in type $Q(n-1)$). Finally, we take $\fe_{\1}=\fc_{\1}$. 

\addtocounter{footnote}{0}\footnotetext{For $\iz:=(n+1)/2$, the negative of the root $\b_{\iz}=\ep_i+\ep_{n+1-i}=2\ep_{\iz}$ is not a root in $P(n-1)$, so terms involving $-\b_{\iz}$ in the definitions of $x_0$ and $\fe_1$ are to be ignored. For example, $x_0=\sum_{i=1}^{\iz-1} d_i(x_{\b_i}+x_{-\b_i}) + d_{\iz} x_{\b_{\iz}}$.}
%\stepcounter{footnote}\footnotetext{See text for explanation.}

\subsection{Dadok-Kac Calculations} \Label{SS:DK}
Except for the simple Lie superalgebras of types $A(n,n)$ and $P(n)$, one can check (see representative example, below) that 
\begin{equation} \Label{E:c1incx0}
\fc_{\1} \subseteq \fc_{x_0}:=\left\{\, y \in \fg_{\1}\mid [\fg_{\0},\, y] \subseteq [\fg_{\0},\, x_0] \,\right\}
\end{equation}
(recall \eqref{E:cv}). Moreover, using Table \ref{T:cohomrings}, $\dim \fc_{\1}=\dim S(\fg_{\1}^*)^{G_{\0}}$, and so it follows from the discussion in Section \ref{SS:dadokkac} that $\fc_{\1}=\fc_{x_0}$ is a Cartan subspace and the action of $G_{\0}$ on $\fg_{\1}$ is polar.

For example, consider the superalgebra $\gl$. Since $\mathfrak{gl}(m|n) \cong \mathfrak{gl}(n|m),$ one may assume $m\le n$. Then $\Omega=\setof{\b_i:=\ep_i-\ep_{m+i}}{1\le i\le m}$ and we write $x_0=\sum_{i=1}^m d_i (x_{\b_i}+x_{-\b_i})$ (with $d_i\ne 0$ for all $i$ and $d_i^2\ne d_j^2$ for $i\ne j$). Conditions (1)--(3) of Theorem \ref{T:closedorbits} are immediate for the set of weights $\Omega\cup -\Omega$ appearing in $x_0$. To check \eqref{E:c1incx0} it suffices to verify that for every root vector $z=x_\ga\in\fg_{\0}$ or $z\in\fh$, and for every $y\in\fc_{\1}$, there exists $z'\in \fg_{\0}$ such that \begin{equation} \Label{E:zyzx0}
[z,y]=[z',x_0].
\end{equation}
In fact, because of the symmetry of $y$ and $x_0$ with respect to positive and negative root vectors, we need only consider $\ga\in\Phi_{\0}^+$. Write $y=\sum_{i=1}^m c_i (x_{\b_i}+x_{-\b_i})$.

Suppose first that $\ga=\ep_j-\ep_k$ with $1\le j<k \le m$, so that $z=E_{j,k}$ (the standard matrix unit). Consider $z'=aE_{j,k}+bE_{m+j,m+k}$ with $a, b\in \C$ to be determined. Then \eqref{E:zyzx0} becomes
\[
\left(
\begin{array}{cc}
d_k  &  -d_j   \\
-d_j  &  d_k   
\end{array}
\right)
\left(
\begin{array}{c}
a \\
b
\end{array}
\right)
=
\left(
\begin{array}{c}
c_k \\
-c_j
\end{array}
\right),
\]
which has a solution since $d_k^2-d_j^2 \ne 0$. A similar analysis applies to $\ga=\ep_{m+j}-\ep_{m+k}$.

Next, if $\ga=\ep_{m+j}-\ep_k$ with $1\le j \le m$ and $2m<k\le m+n$, so that $z=E_{m+j,k}$, then $z'=-(c_j/d_j)z$ satisfies \eqref{E:zyzx0}.  And if $\ga=\ep_j-\ep_k$ with $2m<j<k\le m+n$ then $[z,y]=0$.

Finally, if $z=E_{j,j}\in\fh$ ($1\le j\le m+n$), then taking $z'=(c_j/d_j)z$ if $1\le j\le m$, $z'=(c_{j-m}/d_{j-m})z$ if $m<j\le 2m$, and $z'=0$ if $2m<j\le m+n$ satisfies \eqref{E:zyzx0}.

For the corresponding simple Lie superalgebra $\slmn$ with $m<n$, the same definitions of $\Omega, \fc_1$, and $x_0$ work, and the only changes to the proof involve the situation $z\in\fh$, due to the condition $\str(z)=\str(z')=0$. But we can take advantage of the fact that $[E_{k,k},x_0]=0$ for $2m<k\le m+n$ to find, for any $z\in\fh$, a $z'\in\fh$ satisfying \eqref{E:zyzx0}.  The remaining details are left to the reader.

However, for $\mathfrak{sl}(n|n)$ this argument breaks down, and in fact one can show that there is no Cartan subspace, so the action is not polar. A similar situation arises with the simple Lie superalgebra of type $P(n-1)$.  In both cases the action is stable.  If it were also polar, then the Cartan subspace would contain a generic element $x_{0}.$  By maximality of orbit dimension, one sees that $\c_{x_{0}}$ would equal the Cartan subspace.  However, by a direct calculation we verified that in these two cases $\dim \c_{x_{0}}=1 < \dim S(\g_{\1}^{*})^{G_{\0}},$ so the action cannot be polar.  For these two Lie superalgebras we use the Luna-Richardson theory instead.
\medskip

\renewcommand{\arraystretch}{1.2}
\begin{table}[htdp]
\caption{Dimensions of classical Lie superalgebras}
%\vskip 2in
%\rotatebox{90}{
%\vbox{
\begin{center}
\begin{tabular}{ccc}
$\fg$ &  $\dim \fg_{\0}$ & $\dim \fg_{\1}$   \\ \hline
$\gl$ & $m^2+n^2$ & $2mn$   \\
$\slmn$, $m\ne n$& $m^2+n^2-1$ & $2mn$ \\
%\raisebox{.5ex}{$m\ne n$} \\
$\slnn$ &  $2n^2-2$ & $2n^2$  \\
$\mathfrak{osp}(2m+1|2n)$ & $2m^2+m+2n^2+n$ & $(2m+1)(2n)$  \\
$\mathfrak{osp}(2m|2n)$ & $2m^2-m+2n^2+n$ & $4mn$  \\
$P(n-1)$ & $n^2-1$ & $n^2$  \\
$\widehat Q(n-1)$ & $n^2$ & $n^2$  \\
$Q(n-1)$ & $n^2-1$ & $n^2-1$   \\
$D(2,1;\a)$ & 9 & 8  \\
$G(3)$ & 17 & 14  \\
$F(4)$ & 24 & 16  
\end{tabular}
\end{center}
%}
%}
\Label{T:dimensions}
\end{table}%

\renewcommand{\arraystretch}{1.2}
\begin{table}[htdp]
\caption{Centralizer of $x_0$}
%\vskip 2in
%\rotatebox{90}{
%\vbox{
\begin{center}
\begin{tabular}{ccc}
$\fg$  & $H$ & $\dim H$  \\ \hline
$\gl$  & $T^r \times GL_{|n-m|}$ & $r+(n-m)^2$  \\
$\slmn$, $m\ne n$&  $T^r \times SL_{|n-m|}$ & $r+(n-m)^2-1$ \\
%\raisebox{.5ex}{$m\ne n$} \\
$\slnn$  & $T^{n-1}$ & $n-1$ \\
$\mathfrak{osp}(2m+1|2n)$, $m\ge n$  & $T^r \times SO_{2(m-n)+1}$ & $r+2(m-n)^2+m-n$ \\
$\mathfrak{osp}(2m+1|2n)$, $m< n$  & $T^r \times Sp_{2(n-m)}$ & $r+2(n-m)^2+n-m$ \\
$\mathfrak{osp}(2m|2n)$, $m>n$  & $T^r \times SO_{2(m-n)}$ & $r+2(n-m)^2+(n-m)$ \\
$\mathfrak{osp}(2m|2n)$, $m\le n$  & $T^r \times Sp_{2(n-m)}$ & $r+2(n-m)^2+(n-m)$ \\
$P(n-1)$   & $T^{\lfloor n/2 \rfloor}$ & $\lfloor n/2 \rfloor$ \\
$\widehat Q(n-1)$  & $T^n$ & $n$  \\
$Q(n-1)$  & $T^{n-1}$  & $n-1$  \\
$D(2,1;\a)$  & $T^2$ & 2  \\
$G(3)$  & $SL_2 \times T$ & 4  \\
$F(4)$  & $SL_3 \times T$ & 9 
\end{tabular}
\end{center}
%}
%}
\Label{T:centralizer}
\end{table}%

\renewcommand{\arraystretch}{1.2}
\begin{table}[htdp]
\caption{Polar and stable classical Lie superalgebras}
%\vskip 2in
%\rotatebox{90}{
%\vbox{
\begin{center}
\begin{tabular}{ccccc}
$\fg$  & Polar & Stable & $W$ & $J$\\ \hline
$\gl$ &  \text{Yes} & \text{Yes} & $\Sigma_r \ltimes \Z_2^r$ & $x_1\dots x_r\prod_{i<j}(x_i^2-x_j^2)$\\
$\slmn$,  $m\ne n$& \text{Yes} & \text{Yes} & $\Sigma_r \ltimes \Z_2^r$ & $x_1\dots x_r\prod_{i<j}(x_i^2-x_j^2)$\\
%\raisebox{.5ex}{$m\ne n$} \\
$\slnn$ & \text{No} & \text{Yes}&---&---\\
$\mathfrak{osp}(2m+1|2n)$, $m\ge n$ & \text{Yes} & \text{Yes} & $\Sigma_r \ltimes \Z_4^r$ & $x_1\dots x_r\prod_{i<j}(x_i^4-x_j^4)$\\
%\raisebox{.5ex}{$m\ge n$} \\
$\mathfrak{osp}(2m+1|2n)$, $m<n$& \text{Yes} & \text{No} & $\Sigma_r \ltimes \Z_4^r$ & $x_1\dots x_r\prod_{i<j}(x_i^4-x_j^4)$\\
%\raisebox{.5ex}{$m<n$} \\
$\mathfrak{osp}(2m|2n)$, $m>n$& \text{Yes} & \text{Yes} & $\Sigma_r \ltimes \Z_4^r$ & $x_1\dots x_r\prod_{i<j}(x_i^4-x_j^4)$\\
%\raisebox{.5ex}{$m>n$} \\
$\mathfrak{osp}(2m|2n)$, $m\le n$& \text{Yes} & \text{Yes} & $\Sigma_r \ltimes (\Z_4^r)_e$ & $x_1^2\dots x_r^2\prod_{i<j}(x_i^4-x_j^4)$\\
%\raisebox{.5ex}{$m\le n$} \\
$P(n-1)$ & \text{No} & \text{Yes} &---&---\\
$\widehat Q(n-1)$ & \text{Yes} & \text{Yes} & $\Sigma_n$ & $\prod_{i<j} (x_i-x_j)$\\
$Q(n-1)$ & \text{Yes} & \text{Yes} & $\Sigma_n$ & $\prod_{i<j} (x_i-x_j)$\\
$D(2,1;\a)$ & \text{Yes} & \text{Yes} & $Z_4$ & $x_1^3$\\
$G(3)$ & \text{Yes} & \text{Yes} & $Z_4$ & $x_1^3$ \\
$F(4)$ & \text{Yes} & \text{Yes} & $Z_4$ & $x_1^3$ 
\end{tabular}
\end{center}
%}
%}
\Label{T:polarstable}
\end{table}%

Write $H=G_{\0,x_0}$.  Using \cite[Theorem 7]{Sch:78} one sees that $x_0$ is regular provided
$$
\dim H = \dim \fg_{\0} - \dim \fg_{\1} + \dim S(\fg_{\1}^*)^{G_{\0}}.
$$
In Table~\ref{T:dimensions} we give the dimensions of $\fg_{\0}$ and $\fg_{\1}$ for the classical Lie superalgebras of Sections~\ref{SS:typeA}--\ref{SS:typeF4}. Table~\ref{T:centralizer} gives the structure of $H$ and its dimension. (In the table, $T$ denotes the torus $\C^\times$.) In particular, we observe that in every case, $H$ is connected. Using these data, along with the dimensions of the cohomology rings from Table~\ref{T:cohomrings}, one can check that $x_0$ is regular (and thus the action of $G_{\0}$ on $\fg_{\1}$ is stable) in all but one case. Table~\ref{T:polarstable} lists which cases are polar or stable. In the polar cases, we give the structure of the Dadok-Kac group $W=\operatorname{Norm}_{G_{\0}}(\c_{\1})/\operatorname{Stab}_{G_{\0}}(\c_{\1})$ along with its associated Jacobian $J$ (cf.\ Section~\ref{SS:dadokkac}). The notation $(\Z_4^r)_e$ means the subgroup of $r$-tuples of elements of $\Z_4 = \langle \sqrt{-1} \rangle$ having an even number of entries $\pm \sqrt{-1}$. In each case the action of $W$ on coordinates on $\fc_{\1}\simeq \C^r$ is the obvious one, from which the description of the rings of invariant polynomials in Table~\ref{T:cohomrings} can be verified.

\subsection{Luna-Richardson Calculations} \Label{SS:LR}
Here we apply the Luna-Richardson theory to $\fg=\mathfrak{psl}(n|n)$, to compute $S(\fg_{\1}^*)^{G_{\0}}$. We have $G_{\0} \simeq SL(n)\times SL(n)$ and $\fg_{\1} \simeq M_n(\C)\oplus M_n(\C)$, with action $(A,B)\cdot (X,Y)=(AXB^{-1},BYA^{-1})$ ($A, B\in SL(n),\ X, Y\in M_n(\C)$). Let $\Omega=\{\b_1,\dots,\b_n\} \subset\Phi_{\1}^+$ and $x_0\in \fg_{\1}$ be as defined in the previous subsection for $\gl$ with $m=n$. Then $H:=G_{\0,x_{0}}=\{(D,D)\in G_{\0} \mid D \text{ is diagonal} \}$, and $\fe_{\1}:=\fg_{\1}^H = \{ \sum_{j=1}^n u_j x_{\b_j} + v_j x_{-\b_j} \mid u_j, v_j \in \C \}$. An element $(A,B)$ of the normalizer $N:=N_{G_{\0}}(H)$ is determined by a permutation $\s\in\Sigma_n$ and scalars $a_j, b_j\in\C^\times$ with $\prod a_j = \prod b_j = \sgn(\s)$, where $A_{j,\s(j)}=a_j,\ B_{j,\s(j)}=b_j$ for $1\le j \le n$, and all other entries of $A$ and $B$ are 0. Then $(A,B)$ acts on a pair $(u,v)\in(\C^n)^2$ parametrizing an element of $\fe_{\1}$ by sending it to $(u',v')$ where $u'_j=a_j b_j^{-1} u_{\s(j)},\ v'_j=a_j^{-1} b_j v_{\s(j)}$. 

Abusing notation and viewing $u_j, v_j$ as coordinate functions on $\fe_{\1}$, we see that $S(\fe_{\1}^*)^N \simeq \C[f_2,\dots,f_{2n-2},g_n,h_n]$, where $f_2,\dots,f_{2n}$ are the elementary symmetric polynomials in $u_1 v_1,\dots,u_n v_n$, $g_n:=u_1\dots u_n$, and $h_n:=v_1\dots v_n$. Note that $f_{2n}=g_n h_n$ so that $f_{2n}$ is redundant, but this is the only relation.  We conclude, by the Luna-Richardson Theorem, that $S(\fg_{\1}^*)^{G_{\0}}$ is a polynomial algebra in $n+1$ generators in degrees $2, 4, \dots, 2n-2; n, n$. (This result is apparently known to the experts in classical invariant theory \cite{Sch:06}, but we could not find it in the literature. Gruson \cite{Gru:00} has an incorrect description of $S(\fg_{\1}^*)^{G_{\0}}$.  It is based on her earlier computation \cite{Gru:97} in the case $m\ne n$, but seems to not take into account the subtleties of the supertrace zero condition when $m=n$.)

A similar analysis can be carried out for the simple Lie superalgebra of type $P(n-1)$, as has been done by Gruson \cite{Gru:00}. The degrees of the generators of the polynomial invariants on $\fe_{\1}$ are given in Table \ref{T:cohomrings}.  Alternately, the fact that $S(\fg_{\1}^*)^{G_{\0}}$ is a polynomial algebra can be deduced from \cite[Table 1a, lines 16-17]{Sch:78}, and the degrees of its generators can be determined using \cite[Table 1b]{Sch:78} by omitting those generators having positive degree on $\varphi_1$ or $\varphi_1^*$.

\bibliographystyle{amsmath}
%\bibliographystyle{amsalpha}
%\bibliography{BKN}

\def\Dbar{\leavevmode\lower.6ex\hbox to 0pt{\hskip-.23ex \accent"16\hss}D}

\end{document}